\documentclass{article}

\usepackage[T1]{fontenc}
\usepackage[utf8]{inputenc}
\usepackage[english]{babel}
\usepackage{newtxtext}

\usepackage[dvipsnames]{xcolor}

\usepackage{amsmath}
\usepackage{amssymb}
\usepackage{mathtools}
\usepackage{bm}
\usepackage{bbm}
\usepackage{spalign}
\usepackage{enumerate}
\usepackage{fix-cm}
\usepackage{paralist}
\usepackage{adjustbox}

\usepackage{anyfontsize}
\usepackage{exscale,relsize}
\usepackage{tikz}
\usepackage{pgfplots}
\usepgfplotslibrary{fillbetween}
\usepackage{caption}
\usepackage{float}
\usepackage{wrapfig}
\pgfplotsset{minor tick style={draw=none}}
\usepackage{listings}
\usepackage{scalefnt}

\usepackage{charter}
\newcommand{\package}{\texttt{FrankWolfe.jl}\xspace}

\usepackage[obeyspaces]{url}
\usepackage{subcaption}

\usepackage{graphicx}
\usepackage{float}

\graphicspath{{images/}}

\usepackage{fancyvrb}

\usepackage{jlcode}

\usepackage[section]{algorithm}
\usepackage{algorithmicx}
\usepackage{algpseudocode}

\usepackage{xspace}

\usepackage[figuresright]{rotating}
\usepackage{csvsimple}
\usepackage{makecell}
\usepackage{caption}

\definecolor{darkgreen}{rgb}{0,0.5,0}

\usepackage[colorinlistoftodos]{todonotes}
\usepackage[scaled=0.9]{FiraMono}

\usepackage{booktabs}

\usepackage{amsthm}

\newif\ifarxiv
\arxivtrue

\title{Improved algorithms and novel applications\\ of the FrankWolfe.jl library}
\author{\name Mathieu Besançon\thanks{Corresponding authors. All authors are listed alphabetically.} \email \href{mailto:mathieu.besancon@inria.fr}{mathieu.besancon@inria.fr} \\
    \addr Université Grenoble Alpes, Inria, Laboratoire d'Informatique de Grenoble
    \AND
    \name Sébastien Designolle \email \href{mailto:designolle@zib.de}{designolle@zib.de} \\
    \addr Zuse Institute Berlin, Germany
    \AND
    \name Jannis Halbey \email \href{mailto:halbey@zib.de}{halbey@zib.de}\\
    \addr Zuse Institute Berlin, Germany \\
    Technische Universit\"at Berlin, Germany
    \AND
    \name Deborah Hendrych \email \href{mailto:hendrych@zib.de}{hendrych@zib.de}\\
    \addr Zuse Institute Berlin, Germany \\
    Technische Universit\"at Berlin, Germany
    \AND
    \name Dominik Kuzinowicz \email \href{mailto:kuzinowicz@zib.de}{kuzinowicz@zib.de}\\
    \addr Zuse Institute Berlin, Germany \\
    Technische Universit\"at Berlin, Germany
    \AND
    \name Sebastian Pokutta\footnotemark[1] \email \href{mailto:pokutta@zib.de}{pokutta@zib.de} \\
    \addr Zuse Institute Berlin, Germany \\
    Technische Universit\"at Berlin, Germany
    \AND
    \name Hannah Troppens \email \href{mailto:Troppens@zib.de}{troppens@zib.de}\\
    \addr Zuse Institute Berlin, Germany
    \AND
    \name Daniel Viladrich Herrmannsdoerfer \email \href{mailto:viladrich@zib.de}{viladrich@zib.de} \\
    \addr Zuse Institute Berlin, Germany
    \AND
    \name Elias Wirth \email \href{mailto:wirth@math.tu-berlin.de}{wirth@math.tu-berlin.de}\\
    \addr Technische Universit\"at Berlin, Germany\\
    \addr Zuse Institute Berlin, Germany
}

\DeclareMathOperator*{\argmax}{argmax}
\DeclareMathOperator*{\argmin}{argmin}
\DeclareMathOperator*{\tr}{tr}
\DeclareMathOperator*{\diag}{diag}
\newcommand{\innp}[2]{\left\langle #1, #2 \right\rangle}
\newcommand{\norm}[1]{\left\| #1 \right\|}
\newcommand{\vx}{\mathbf{x}}
\newcommand{\vs}{\mathbf{s}}
\newcommand{\vw}{\mathbf{w}}
\newcommand{\vvv}{\mathbf{v}}

\newcommand{\vy}{\mathbf{y}}
\newcommand{\vd}{\mathbf{d}}

\newcommand{\vl}{\mathbf{l}}
\newcommand{\vu}{\mathbf{u}}
\newcommand{\va}{\mathbf{a}}
\newcommand{\vb}{\mathbf{b}}
\newcommand{\vg}{\mathbf{g}}

\renewcommand{\leq}{\leqslant}

\renewcommand{\geq}{\geqslant}

\newcommand{\Xset}{{\ensuremath{\mathcal{X}}}\xspace}
\newcommand{\Fcal}[1]{{\ensuremath{\mathcal{F}(#1)}}\xspace}

\usepackage{jmlr2e}
\usepackage[noabbrev,capitalize,nameinlink]{cleveref}[0.21]
\usepackage{natbib}
\setcitestyle{authoryear,round,citesep={;},aysep={,},yysep={;}}
\usepackage[vvarbb]{newtxmath}
\usepackage{multirow}

\definecolor{labelkey}{rgb}{0,0.08,0.45}
\definecolor{refkey}{rgb}{0,0.6,0.0}
\definecolor{Brown}{rgb}{0.45,0.0,0.05}
\definecolor{dgreen}{rgb}{0.00,0.49,0.00}
\definecolor{dblue}{rgb}{0,0.08,0.75}
\definecolor{ffwwqq}{rgb}{1.,0.4,0.}
\definecolor{qqzzqq}{rgb}{0.,0.6,0.}
\definecolor{qqqqff}{rgb}{0.,0.,1.}
\definecolor{dred}{HTML}{D90404}
\definecolor{orng}{HTML}{D35400}
\definecolor{cb-black}      {RGB}{  0,   0,   0}
\definecolor{cb-blue-green} {RGB}{  0,  073,  073}
\definecolor{cb-green-sea}  {RGB}{  0, 146, 146}
\definecolor{cb-rose}       {RGB}{255, 109, 182}
\definecolor{cb-salmon-pink}{RGB}{255, 182, 119}
\definecolor{cb-purple}     {RGB}{ 73,   0, 146}
\definecolor{cb-blue}       {RGB}{ 0, 109, 219}
\definecolor{cb-lilac}      {RGB}{182, 109, 255}
\definecolor{cb-blue-sky}   {RGB}{109, 182, 255}
\definecolor{cb-blue-light} {RGB}{182, 219, 255}
\definecolor{cb-burgundy}   {RGB}{146,   0,   0}
\definecolor{cb-brown}      {RGB}{146,  73,   0}
\definecolor{cb-clay}       {RGB}{219, 209,   0}
\definecolor{cb-green-lime} {RGB}{ 36, 255,  36}
\definecolor{cb-yellow}     {RGB}{255, 255, 109}
\definecolor{bred}{HTML}{FF0000}
\definecolor{bpurp}{HTML}{BF00BF}
\definecolor{bblu}{HTML}{0000FF}
\definecolor{bcyan}{HTML}{00BFBF}
\definecolor{byellow}{HTML}{BFBF00}
\definecolor{bgreen}{HTML}{008000}

\begin{document}

\maketitle

\begin{abstract}
Frank-Wolfe (FW) algorithms have emerged as an essential class of methods for constrained optimization, especially on large-scale problems.
In this paper, we summarize the algorithmic design choices and progress made in the last years of the development of FrankWolfe.jl, a Julia package gathering high-performance implementations of state-of-the-art FW variants.
We review key use cases of the library in the recent literature, which match its original dual purpose: first, becoming the de-facto toolbox for practitioners applying FW methods to their problem, and second, offering a modular ecosystem to algorithm designers who experiment with their own variants and implementations of algorithmic blocks.
Finally, we demonstrate the performance of several FW variants on important problem classes in several experiments, which we curated in a separate repository for continuous benchmarking.
\end{abstract}

\section{Introduction}

In the last decade, Frank-Wolfe (FW) algorithms have gained significant traction in various applications at the intersection of constrained and nonlinear optimization.
We present a summary of progress made over the last years of the \package{} package, including recent advances in FW algorithms that have been integrated into the package,
its extensions developed as companion packages, and selected applications from the literature that leveraged \package{} for computational experiments or further development of other software libraries.\\

\package{} was initially introduced in \citet{besanccon2022frankwolfe} after its release in 2021 and we refer interested readers to this article for prominent features of the first releases.
In this paper, we will primarily highlight progress made after this initial publication and available in the 0.4.7 version.
The rest of this paper is structured as follows.
\cref{sec:fwalgs} introduces the new FW-based algorithms implemented in the package, \cref{sec:stepsize} reviews the different step-size strategies introduced from the literature,
\cref{sec:lmo} details new feasible regions and corresponding linear minimization oracles.
\cref{sec:quadratic} presents the new active set abstraction and different implementations for specific use cases.
\cref{sec:dynamic} introduces mechanisms to dynamically modify the algorithm run.
The appendix contains the pseudo-code for the main new FW variants, a selection of recent applications of the package in the literature, and a benchmark suite and results we developed to monitor the performance of the different algorithms.

Throughout the paper and unless specified otherwise, we will consider optimization problems of the form:
\begin{align*}
    \min_{\vx \in \Xset{}}\; & f(\vx),
\end{align*}
with $\Xset{}$ a closed, convex subset of a Hilbert space and $f: \Xset{} \rightarrow \mathbb{R}$ a differentiable objective function.
The iteration count of algorithms will be denoted with $t$ and the iterate at the current iteration as $\vx_t$.
For $n \in \mathbb{N}$, we also note $[n] := \{1, \dots, n\}$. We denote the unit simplex in $\mathbb{R}^n$ as $\Delta_n$.
We will now briefly outline the algorithmic improvements and developments in \package{}, including new algorithms, step-size strategies, and oracles.

\section{New Frank-Wolfe algorithms}\label{sec:fwalgs}

In the following, we provide an overview of new FW variants implemented in \package{}.
For details on the algorithms and theoretical guarantees, we refer the interested reader to the original papers.
The pseudo-code for the different algorithms are presented in \cref{sec:pseudocode}.

\subsection{Blended Pairwise Conditional Gradients (BPCG)}

The \emph{Blended Pairwise Conditional Gradients} (BPCG) algorithm from \citet{tsuji2022pairwise} is an improved variant of the \emph{Pairwise Conditional Gradients} (PCG) method.
Unlike PCG, BPCG eliminates swap steps, which may hinder primal progress and limit applicability in infinite-dimensional settings.
BPCG combines PCG with the blending criterion from the Blended Conditional Gradients (BCG) algorithm \citep{braun2019blended}, leading to better convergence rates and improved performance.
It is straightforward to implement, avoids complex gradient alignment procedures (such as the ones, e.g., in \citet{rinaldi2020unifying,rinaldi2023avoiding}), and ensures solutions are sparse, making it particularly effective for applications like sparse signal recovery and kernel herding.
Numerical experiments in \citet{tsuji2022pairwise} demonstrate that BPCG produces sparser solutions than PCG, offering significant advantages in computational efficiency and performance, which is in line with the performance evaluation of our implementation in \package{} and of the benchmarks presented in \cref{sec:benchmarks}.

The pseudo-code is provided in Algorithm~\ref{alg:BPCG} in \cref{sec:pseudocode}. The implementation in \package{} uses the relaxed criterion with $\kappa = 2.0$ from \cite[Section 3.4]{tsuji2022pairwise} in Line~\eqref{bpcg:select} instead of $\kappa = 1.0$, for significantly enhanced sparsity while maintaining almost identical convergence speed in iteration and time.
In the description of Algorithm~\ref{alg:BPCG} we use $c[\vx_{t}](\va_{t})$ to denote the weight of $\va_t$ in the convex decomposition of $\vx_t$ that the algorithm maintains.
This is the current default FW variant if an active set can be maintained.
We highlight that sparsity of the active set provides multiple algorithmic advantages, including faster convergence (see \cref{sec:pivoting}) and lower iteration costs, since active set operations usually have a cost linear in the number of vertices.

\subsection{Pairwise Conditional Gradients (PCG)}

We added an implementation of the \emph{Pairwise Conditional Gradients} (PCG) from \citet{lacoste15} similar to the Away-step Frank-Wolfe algorithm in \package{}.
PCG performs a pairwise step which consists in a direction $\vd_t := \va_t - \vvv_t$ following notation of \cref{alg:BPCG}.
The algorithm also features an optional lazification mechanism introduced in \citet{braun2019lazifying} similar to the other ones, i.e., the possibility to avoid expensive calls to the exact linear minimization oracle (LMO) at every iteration, using vertices readily available for directions and maintaining an estimate of the dual gap.

\subsection{Decomposition-Invariant Conditional Gradients (DICG)}

The \emph{Decomposition-Invariant Conditional Gradients} (DICG) corresponds to the algorithm defined in \citet{garber2016linear} for the case of ``simplex-like'' polytopes and \citet{bashiri2017decomposition} for general polytopes.
The motivation for the algorithm is two-fold.
First, typically fast FW variants on polytopes rely on an active set which induces a storage cost at least proportional to the active set complexity of the optimal solution and to the size of the in-memory representation of the vertices.
Second, the evolution of the convex decomposition of the iterates throughout the algorithm depends on the steps that are performed and the set of directions available for local steps (i.e., steps depending on the active set) depends on the current convex decomposition.
DICG thus aims at building decomposition-invariant iterations that do not rely on storage nor exploitation of a set of active vertices.
We consider the case of polytopes inscribed in the 0-1 hypercube for simplicity of exposure.
We will denote the current face of the iterate $\vx$:
\begin{align*}
    \Fcal{\vx} = \{\vy \in \Xset, y_i = x_i~\forall i \in [n], x_i \in \{0,1\}\}.
\end{align*}

The key aspect of DICG is considering ``extended'' linear minimization oracles that can compute extreme points of a face of the polytope.
In many cases, such a restricted oracle is not computationally more expensive than the original oracle and will reduce the algorithm complexity compared to the original LMO.
As an example, the standard LMO for the Birkhoff polytope consists of performing the Hungarian algorithm with the costs provided by the gradient.
The decomposition-invariant oracle also requires being able to fix a variable to zero, which consists in forbidding a row/column pair in the Hungarian algorithm,
or to fix a variable to one, which removes a row and a column and performs the algorithm on a submatrix.

The initial convergence analysis of DICG from \citet{garber2016linear} was complemented in \citet{wirth2024fast} with an affine-independent convergence rate under an assumption of error bound on the objective function.

\subsection{Blended Decomposition-Invariant Conditional Gradients}

The \emph{Blended Decomposition-Invariant Conditional Gradients} (BDICG) algorithm is the natural combination of the DICG algorithm with the blending
used in the BPCG algorithm, resulting in an even faster DICG variant and bearing resemblance to the in-face Frank-Wolfe algorithm of \citet{freund2017extended},
however without relying on expensive oracles. The pseudo-code is presented in \cref{alg:bdicg}.
Furthermore, decomposition-invariant algorithms naturally fit optimization problems on polytopes, while it is not obvious they would provide much or any benefit on other convex sets which were the focus of the in-face Frank-Wolfe algorithm.

The convergence of BDICG has been studied in an affine-invariant setting in \citet{wirth2024fast} under a sharpness assumption of the objective function on polytopes, the established rate is sublinear, improves over the classic $\mathcal{O}(t^{-1})$ in the general case, and recovers the linear rate in the case of strongly convex functions that is common for FW variants with enhanced directions (away or blended pairwise steps).

\subsection{Block-Coordinate Frank-Wolfe (BCFW)}

The \emph{Block-Coordinate Frank-Wolfe} (BCFW) algorithm is designed for optimization problems with block-separable constraints.
Our implementation allows for a custom update order unifying the cyclic and random block variants proposed in \citet{beck2015cyclic} and \citet{lacostejulien2013blockcoordinate}.
Additionally, the implementation offers flexibility in the updating mechanism of different blocks, including the update direction and step size.
In addition to the standard FW step used in the original papers, we also provide a blended pairwise step.

Another Frank-Wolfe-based algorithm with a block structure was proposed in \citet{woodstock2024flexible} and its convergence rate established in the convex and nonconvex settings.
In contrast with previous work on FW with a block structure, it considers partial block updates under a cyclic assumption from the proximity operator literature.
This notably allows for provable convergence while updating blocks at different frequencies, e.g., based on the computational cost of the corresponding LMOs.
The numerical experiments in \citet{woodstock2024flexible} were conducted with the new algorithm implemented on top of \package{}.

\begin{algorithm}
    \caption{Block-Coordinate Frank-Wolfe (BCFW)}\label{alg:bcfw}
    \begin{algorithmic}[1]
        \Require Convex smooth function $f \colon \bigtimes_{i=1}^m \mathcal{X}_i \to \mathbb{R}$, starting vertex $\vx_0 \in \bigtimes_{i=1}^m \mathcal{X}_i$, step size $\gamma$.
        \Ensure Points $\vx_{1}, \ldots, \vx_T$ in $\bigtimes_{i=1}^m \mathcal{X}_i$.
        \For{$t \in 0\dots T$}
            \State Choose set of blocks $I_t \subset \{1,\dots,m\}$.
            \State $\vg_t \gets \nabla f(\vx_t)$
            \For{$i=1$ \textbf{to} $m$}
            \If{$i \in I_t$}
                \State $\vvv_t^i \gets \argmin_{\vvv \in \mathcal{X}_i} \innp{\vg_t^i}{\vvv}$
                \State $\vx_{t+1}^i \gets \vx_t^i + \gamma_{t}^i (\vvv_t^i - \vx_t^i)$
            \Else
                \State $\vx_{t+1}^i \gets \vx_t^i$
            \EndIf
            \EndFor
        \EndFor
        \State \Return $\vx_T$
    \end{algorithmic}
\end{algorithm}

\subsection{Alternating Linear Minimization (ALM)}

The \emph{Alternating Linear Minimization} (ALM) algorithm is inspired by von Neumann's alternating projections algorithm for solving feasibility problems over intersections of convex sets, see \citet{braun2023alternating}. However, unlike von Neumann's original algorithm, the algorithm utilizes linear minimization oracles to avoid projections.
The key idea of ALM lies in rewriting the feasibility problem as an optimization problem over the product of the given sets. This reformulation also allows us to solve optimization problems over intersections of convex sets by simply extending the objective function.
Therefore, ALM enables combining different constraints and their corresponding LMOs in a single optimization problem.
Our implementation is up to minor modifications identical to the Split Conditional Gradients (SCG) algorithm from \citet{woodstock2024splitting}.

\subsection{A note on dual gaps}

Many methods for lazyfying textbook algorithms use an upper bound $\Phi$ on the \emph{strong FW gap} rather than the standard \emph{FW gap} as a criterion to select the desired step.
However, throughout the package, we estimate the standard FW gap via $\Phi$ and not the strong FW gap.
The reason is to be able to use $\Phi$ as a stopping criterion in a fashion consistent with the standard FW gap; moreover strong convexity parameters (and pyramidal width parameters) required for geometric strong convexity bounds relying on the strong FW gap are usually unknown.
Our choice in the package is justified by the argumentation in \citet[Proof of Theorem 3.1]{braun2019blended}, which shows that a bound on the standard FW gap induces a bound on the strong FW gap, when algorithms use any steps that imply away steps.\\
With the notation of Algorithms~\ref{alg:BPCG} and \ref{alg:bdicg}: if the FW gap satisfies
\begin{align*}
    \innp{\nabla f(\vx_t)}{\vx_t - \vvv_t} \leq \Phi,
\end{align*}
then the strong FW gap satisfies:
\begin{align*}
    \innp{\nabla f(\vx_t)}{\va_t - \vvv_t} \leq k \Phi,
\end{align*}
with $k$ being a small positive constant, typically $k = 2$.

\section{New step-size strategies}\label{sec:stepsize}

We next present new step sizes from the recent literature which have been implemented in \package{}.
These step sizes have been shown to perform well under particular settings or classes of problems (e.g., polytope constraints, generalized self-concordant functions, or optimum in the interior).

\subsection{Open-loop step sizes}

The FW algorithm with line-search or short-step is known to converge at rates not faster than $\Omega(t^{-1-\epsilon})$
for any $\epsilon>0$ in the setting of \emph{Wolfe's lower bound} \citep{wolfe1970convergence,canon1968tight}.
This setting assumes that the feasible region is a polytope, the optimal solutions lie in the relative interior of an at least one-dimensional face, and the objective is strongly convex.
Recently, open-loop step sizes of the form
\begin{align*}
    \gamma_t & = \frac{\ell}{t+\ell} \qquad \text{for} \ \ell \in\mathbb{N}_{\geq 1}
\end{align*}
have been identified as a means to achieve faster convergence rates of up to $\mathcal{O}(t^{-2})$ in Wolfe's setting \citep{bach2021effectiveness,wirth2023acceleration,wirth2023accelerated}.
This line of research culminated in the development of automatically adapting step-size rules, such as
\begin{align}
    \gamma_t & = \frac{2 + \log(t+1)}{t+2 + \log(t+1)},
\end{align}
that adapt to the geometry of the problem setting to achieve rates most often optimal among open-loop step sizes up to an initial burn-in phase of potential lower convergence speed \citep{pokutta2024frank}.
To facilitate further experimentation, the \package{} package implements all of the above step sizes using the \texttt{Agnostic} and \texttt{GeneralizedAgnostic} step-size options.
It also supports general functions $g\colon \mathbb{N} \to \mathbb{R}_{\geq 0}$, allowing for step sizes of the form
\begin{align}
    \gamma_t & = \frac{g(t)}{t+g(t)}
\end{align}
to facilitate further experimentations.

\subsection{Monotonic open-loop step sizes}

In \citet{carderera2021simple,carderera2021simple2}, the study of FW algorithms for generalized self-concordant functions motivated the introduction of
a modification of the open-loop step size ensuring monotonic primal progress.
The step size consists of:
\begin{align*}
    \gamma_t = \max_{N \in \mathbb{N}_{\geq N_t}}& \gamma(N) := \frac{2^{-N}}{2+t} \\
    \text{s.t. } & \vx_t - \gamma \vd_t \in \mathrm{dom}(f)\\
    & f(\vx_t - \gamma \vd_t) < f(\vx_t),
\end{align*}
where $N_t$ is the number of halving steps performed across previous iterations, starting from $N_1 = 1$.
The corresponding procedure halves the step derived from the open-loop step size until the iterate lies in the domain of $f$ and yields primal progress.
The domain membership criterion may appear redundant with monotonicity. In practice, however, a domain oracle can be called first to short-circuit the evaluation and avoid calls to the objective function, since it is in most cases less demanding than function evaluations.
Despite its simplicity, the monotonic open-loop step size is sufficient to ensure convergence at the ``usual'' rates ($1/t$ in general and linear convergence for e.g. uniformly convex feasible regions or strongly convex functions over polytopes) for generalized self-concordant functions despite the absence of a global Lipschitz constant.

\subsection{Adaptive step sizes}\label{subsubsec:Adaptive}

In \citet{pedregosa2018step}, an adaptive step-size strategy was proposed that dynamically estimates the smoothness constant.
While extremely powerful with great real-world performance, this step-size strategy often suffers from numerical instabilities due to the test criterion employed, which basically relies on repeatedly testing the smoothness inequality for different estimations of the smoothness constant.
In order to overcome these numerical instabilities, an alternative adaptive step-size strategy was proposed in \citet{pokutta2024frank}, which is similar in nature to the one in \citet{pedregosa2018step}, however uses the alternative test:
\begin{equation*}
    \label{eq:altAdaptive}
    \innp{\nabla f(\vx_{t+1})}{\vx_t - \vvv_t} \geq 0,
\end{equation*}
to accept an estimate $\tilde L$, where $\vx_{t+1}$ is chosen according to the short step with smoothness estimate $\tilde L$.
In \citet{pokutta2024frank}, it was shown that this test is equivalent to the original smoothness inequality test from \citet{pedregosa2018step} but provides significantly improved numerical stability.
The adaptive step-size procedure is given in \cref{alg:AdaptiveStepSize}.
This step-size strategy is the default strategy in \package{}.

\begin{algorithm}
    \caption{Modified adaptive step-size strategy}\label{alg:AdaptiveStepSize}
    \begin{algorithmic}[1]
        \Require Objective function $f$, smoothness estimate $\widetilde{L}$, feasible point $\vx$ and direction $\vd$ with $\innp{\nabla f(\vx)}{\vd} \geq 0$, progress parameters $\eta \leq 1 < \tau$, maximum step $\gamma^{\max}$
        \State $M \gets \eta \widetilde{L}$
        \While{true}
        \State $\gamma \gets \min \left\{\frac{\innp{\nabla f(\vx)}{\vd}}{M \norm{\vd}^{2}}, \gamma^{\max} \right\}$ \Comment{compute short step with estimate $M$}
        \If{$\innp{\nabla f(\vx - \gamma \vd)}{\vd} \geq 0$} \Comment{sufficient decrease}
        \State $\widetilde{L}^{*} \gets M$
        \State \Return $\widetilde{L}^{*}$, $\gamma$
        \EndIf
        \State $M \gets \tau M$
        \EndWhile
    \end{algorithmic}
\end{algorithm}

\subsection{Secant step sizes}

We implemented a secant-based line search strategy introduced in \citet{hendrych2025secant} and which works particularly well for generalized self-concordant functions.
Although it requires gradient evaluations during the line search procedure, experimental results showed its promising performance in several scenarios due to its favorable convergence properties (assuming suitable function properties).
We highlight that we added a domain oracle with an initial bound on the maximum step size to enable the use of the line search for generalized self-concordant functions when $\Xset \not\subset \mathrm{dom}(f)$, i.e., one has to verify that a point is within the domain before being able to evaluate a gradient.

\begin{algorithm}
    \caption{Secant step-size strategy}\label{alg:secant}
    \begin{algorithmic}[1]
        \Require Objective function $f$, feasible points $\vx$, direction $\vd$ with $\innp{\nabla f(\vx)}{\vd} \geq 0$, $\gamma^{\max}$
        \While{ $\vx - \gamma^{\max} \vd \notin \mathrm{dom}(f)$ }
        \State $\gamma^{\max} \gets \gamma^{\max} / 2$
        \EndWhile
        \State $\gamma_{-1} \gets 0$, $\gamma_0 \gets \gamma^{\max}$
        \State $\vx_{-1} \gets \vx$, $\vx_0 \gets \vx - \gamma_0 \vd$
        \For{$k \in 1\dots K$}
            \If{$\innp{\nabla f(\vx_{k-1}) - \nabla f(\vx_{k-2})}{\vd} \approx 0$}
                \State \Return $\gamma_{k-1}$
            \EndIf
            \State $\tilde{\gamma}_k \gets \gamma_{k-1} - \innp{\nabla f(\vx_{k-1})}{\vd} \frac{\gamma_k - \gamma_{k-1}}{\innp{\nabla f(\vx_{k-1}) - \nabla f(\vx_{k-2})}{\vd}} $
            \State ${\gamma}_k \gets \mathrm{project}_{\left[0,\gamma^{\max}\right]}(\tilde{\gamma}_k)$
            \State $\vx_k \gets \vx - \gamma_k \vd$
        \EndFor
        \State \Return $(\gamma_{k_{\min}})_{k\in1\dots K}$ minimizing the function value.
    \end{algorithmic}
\end{algorithm}

\section{Linear Minimization Oracles}\label{sec:lmo}

Since the first version of the package, multiple linear minimization oracles were added to support solving different types of problems.

\paragraph*{Spectraplex}

The spectraplex (respectively, unit spectraplex) is the intersection of the cone of positive semidefinite matrices with an equality (respectively, inequality) constraint on the trace of the matrix:
\begin{align*}
    \mathcal{S}_{\tau} = \left\{X \in \mathbb{S}_+^n, \mathrm{Tr}(X) = \tau \right\}.
\end{align*}
The spectraplex is the convex hull of symmetric outer products $\vvv \vvv^\top$ such that $\norm{v} = \sqrt{\tau}$, and given a direction, a minimizing extreme point is given by computing the leading eigenvector of the negative direction.

\paragraph*{Hypersimplex}

The hypersimplex of radius $\tau$ is the convex hull of nonnegative vectors with $K$ coordinates taking value $\tau$ and all others set to zero, and can be expressed as:
\begin{align*}
    \mathcal{H}_{K, \tau} = \left\{\vx \in \left[0,\tau\right]^n, \sum_i x_i = K \tau \right\}.
\end{align*}
The linear oracle consists in setting the $K$ entries corresponding to the smallest gradient values to $\tau$.
A hypersimplex oracle was used in \citet{macdonald2022interpretable} to compute the $K$ most relevant pixels for the classification of images.

\paragraph{Symmetric oracle and symmetric Frank-Wolfe algorithm}

Suppose that there is a group $G$ acting on the underlying vector space containing $\mathcal{X}$ and such that for all $\vx \in\mathcal{X}$ and $g\in G$:
\begin{align*}
    f(g\cdot \vx)=f(\vx)\quad\text{and}\quad g\cdot \vx\in\mathcal{X}.
\end{align*}
Then, the computations can be performed in the subspace invariant under the action of $G$.
Let us indeed consider an optimal solution $\vx^*\in\mathcal{X}$ and its group average (also called Reynolds operator) defined by
\begin{align*}
    \mathcal{R}(\vx)=\frac{1}{|G|}\sum_{g\in G}g\cdot \vx.
\end{align*}
By convexity of $f$ and $\mathcal{X}$, and with the invariance properties above, we have that $\mathcal{R}(\vx^*)\in\mathcal{X}$ and
\begin{align*}
    f\big(\mathcal{R}(\vx^*)\big)=f\bigg(\frac{1}{|G|}\sum_{g\in G}g\cdot \vx^*\bigg)\leq\frac{1}{|G|}\sum_{g\in G}f(g\cdot \vx^*)=f(\vx^*),
\end{align*}
so that $\mathcal{R}(\vx^*)$ is a symmetric optimal solution.

Similarly, the LMO can be performed on $\mathcal{R}(\mathcal{X})$.
If this additional invariance constraint can be exploited to directly accelerate the LMO, one can explicitly work in the invariant subspace.
Otherwise, the package offers a way to pre- and post-process the gradient so that the full LMO can be called within a algorithm working in the reduced subspace.
\begin{jllisting}
    struct SubspaceLMO{LMO,TD,TI} <: LinearMinimizationOracle
        lmo::LMO
        deflate::TD
        inflate::TI
    end
    struct SubspaceVector{HasMultiplicities,T,DT} <: AbstractVector{T}
        data::DT
        vec::Vector{T}
        mul::Vector{T}
    end
\end{jllisting}
Practically, the iterate lives in a reduced subspace, here the image of the Reynolds operator, and the gradient is first inflated before being passed to the full LMO, from which the output is then deflated, which amounts to projecting onto the symmetric subspace in this case.
The speedups obtained have various origins:
\begin{itemize}
    \item The number of iterations is reduced, roughly by a factor corresponding to the average size of the orbits of the outputs of the nonsymmetrized LMO.
    \item For active set-based methods, the active set also benefits from this reduction, which can have a noticeable impact on the speed of the search over it when its size becomes quite large.
    \item When the memory size of the symmetrized iterate is smaller than the original one, basic in-place operations performed at each iteration are also faster, which can make a significant difference for large numbers of lazy iterations.
        For this purpose, the structure \texttt{SubspaceVector} simplifies the physical reduction of the memory footprint of the iterate; see the corresponding example in the documentation, named \emph{Accelerations for quadratic functions and symmetric problems}.
\end{itemize}

Note that this dimension reduction does not only work for symmetric instances, but, generally, for all cases where the space of a \emph{natural} LMO is larger than the subspace in which the computations are effectively happening.
The collaborative filtering below is such a nonsymmetric example, although the benchmarks do not exploit the dimension reduction.

\section{Abstract and specialized active set structures}\label{sec:quadratic}

Many FW variants operate on the explicit decomposition of the current iterate into a convex combination of a subset of vertices $\mathcal{S}$:
\begin{align*}
    (\lambda, \mathcal{S}): \vx = \sum_{k \in \mathcal{S}} \lambda_k \vvv_k.
\end{align*}
Although the active set typically refers to the subset $\mathcal{S}$ with nonzero weights, we use it within \package{} to refer to the tuple $(\lambda, \mathcal{S})$ captured together with $\vx$ in a data structure.
Keeping the three elements in a single structure allows, for instance, for fast updates when few weights are changed, or specialized operations for particular atom types.

A default \texttt{ActiveSet} structure is used in typical calls to the active set-based algorithms.
The active set constitutes the essential part of the state of these algorithms, and they can therefore all be called with an already created active set as argument, for instance, to warm-start them.

\paragraph*{ActiveSetQuadraticProductCaching.}

We also implemented a second active set type for convex quadratic minimization with the objective written as:
\begin{align*}
    f(\vx) = \frac{1}{2} \vx^\top A \mathbf{x} + \mathbf{b}^\top \vx.
\end{align*}
This active set was used in particular in the huge-scale projection problem tackled in~\citet{designolle2023improved}.
The main idea is to avoid the computation of scalar products $\innp{\vvv_1}{\vvv_2}$ or quadratic forms $\innp{\vvv_1}{A \vvv_2}$ within the active set as needed to obtain the away and local vertices $\va_t$ and $\vs_t$ (see, e.g., \cref{alg:BPCG}).
This relies on pre-computing some specific values that can be used, thanks to the linearity of the gradient, to efficiently recover these scalar products at a reduced computational cost.

Formally, if the objective function has a gradient $\nabla f(\vx)=A\vx+\vb$, then, since the current iterate is a convex combination of the atoms of the active set, i.e., $\vx_t=\sum_j c[\vx_t](\vvv_j)\vvv_j$, we have:
\begin{align}
\innp{\nabla f(\vx_t)}{\vvv_i} =\sum_j c[\vx_t](\vvv_j)\innp{A\vvv_j}{\vvv_i} + \innp{\vb}{\vvv_i} &\;\; \forall i \in \mathcal{S}.
  \label{eqn:argminmax}
\end{align}
By storing all values $\innp{A\vvv_j}{\vvv_i}$ and $\innp{\vb}{\vvv_i}$ for $i,j \in \mathcal{S}$, we can maintain an up-to-date value of \cref{eqn:argminmax} according to potential modifications from previous iterations.
Since at most two coefficients $c[\vx_t](\vvv_j)$ can have changed compared to $c[\vx_{t-1}](\vvv_j)$, the cost of this update is small compared to the computation of a scalar product.

\paragraph*{ActiveSetQuadraticLinearSolve.}

In the case of convex quadratic problems, we can further exploit the current active set to gain immediate improvement without additional Frank-Wolfe steps.
At any iteration, we can attempt to minimize $f$ over $\mathrm{conv}(\{\vvv_k\}_{k\in \mathcal{S}})$, which is equivalent to the corrective step performed by the Fully-Corrective Frank-Wolfe (FCFW) algorithm.
Computing this step can however be as hard as the original problem.

Optimality conditions imply that the inner product of the gradient at an optimal point $\vx_\text{opt}$ with any feasible direction has to be nonnegative.
It suffices to check this for vertices in the active set,
$$
    \innp{\nabla f(\vx_\text{opt})}{\vvv_k - \vx_\text{opt}} \geq 0 \qquad \text{for }k \in \mathcal{S}.
$$
Note that this is equivalent to the Frank-Wolfe gap of $f$ over $\mathrm{conv}(\{\vvv_k\}_{k\in \mathcal{S}})$ being nonpositive at $\vx_\text{opt}$.

In \citet{halbey2025efficientquadraticcorrectionsfrankwolfe}, a proposed approach is to relax the problem and minimize over the affine space spanned by the current active vertices $\operatorname{aff}(\{\vvv_k\}_{k\in \mathcal{S}})$.
The gradient has to be orthogonal to the affine space, which leads to equality constraints.
Furthermore, we can use the directions $\vvv_{\tilde k}-\vvv_k$ for a fixed $\tilde k \in \mathcal{S}$ instead of $\vvv_k - \vx_\text{opt}$ as they span the same space,
$$
    \innp{\nabla f(\vx_\text{opt})}{\vvv_k - \vvv_{\tilde k}} = 0 \qquad \text{for }k \in \mathcal{S} \setminus \{\tilde k\}.
$$
Given the convex quadratic objective with gradient $\nabla f(\vx) = A \vx + \vb$ and
using barycentric coordinates $\tilde \lambda$, $\vx_\text{opt} = \sum_{k\in \mathcal{S}} \vvv_k \tilde \lambda_k = V \tilde \lambda$, and
we can rewrite this as:
\begin{align*}
& \innp{A V \tilde \lambda + \vb}{\vvv_k-\vvv_{\tilde k}} = 0 \qquad \text{for }k \in \mathcal{S} \setminus \{\tilde k\}\\
& \sum_{k \in S} \tilde \lambda_k = 1.
\end{align*}
Note that the system consists of $|\mathcal{S}|$ equations with affinely independent columns and has therefore at least one solution.

In order to obtain a point in the convex hull of the active set, we need to recover a nonnegative weight vector $\lambda$.
We can enforce the nonnegativity directly and solve the resulting linear program (LP) with standard solvers even at a large scale, at a cost that is comparable or lower to many LMOs.
In some cases however, this LP is infeasible, namely if the solution space of the linear system does not intersect with the $|\mathcal{S}|$-simplex.

This motivated an alternative direct solving process using \emph{Wolfe's step} used in Wolfe's minimum norm point algorithm \citep{wolfe1976finding}.
One can use the current weights $\lambda'$ and find the point on the segment $[\lambda', \tilde \lambda]$ that is the closest to $\lambda'$ in the simplex through a simple ratio test.
Importantly, this ensures that at least one vertex will be dropped from the active set since the minimum ratio is achieved when at least one entry is zero, assuming $\tilde \lambda$ had negative entries.

This quadratic correction approach is detailed in \citet{halbey2025efficientquadraticcorrectionsfrankwolfe} and implemented as an active set type that can be used on top of either a standard active set or the \texttt{ActiveSetQuadraticProductCaching} type presented before.

\section{Dynamic algorithmic modifications}\label{sec:dynamic}

In this section, we present two mechanisms dynamically modifying the execution of FW algorithms, namely callbacks allowing the user to collect information and stop the algorithm, and the pivoting framework which rewrites iterates.

\subsection{Dynamic stopping criterion with callbacks}

A callback mechanism was introduced in the early versions of \package{} without any algorithmic effect.
This enabled users, e.g., to print additional logs or record information about each iteration.
We added a convention that callbacks can return a Boolean flag indicating whether the algorithm should continue or terminate.
The callback can determine this condition based on the current state, gap, iteration count, time, etc, enabling an adaptive termination.
We illustrate with a use case below in which we terminate BPCG as soon as the active set contains ten atoms.

\begin{lstlisting}[language=julia]
f(x) = 1/2 * norm(x .- 1)^2
grad!(storage, x) = storage .= (x .- 1)
lmo = FrankWolfe.KSparseLMO(5, 2.0)
x0 = FrankWolfe.compute_extreme_point(lmo, ones(100))

function build_callback(max_vertices)
    return function callback(state, active_set, args...)
        return length(active_set.atoms) <= max_vertices
    end
end

# callback initialization with a maximum of ten vertices
callback = build_callback(10)

FrankWolfe.blended_pairwise_conditional_gradient(f, grad!, lmo, x0, callback=callback)
\end{lstlisting}

We also use this example to highlight a common pattern to build callbacks using a closure receiving the callback parameters.
This ensures both flexibility and performance by avoiding the declaration of \texttt{max\_vertices} as a global variable.

\subsection{The pivoting Frank-Wolfe framework}\label{sec:pivoting}

Although FW algorithms empirically maintain a sparse combination of extreme points of the feasible region, this sparsity can rarely be explicitly ensured or constrained.
Indeed, even for a simple least-square objective, an $\ell_1$-norm ball constraint will be sparsity-inducing while keeping the problem tractable, whereas an explicit cardinality constraint would imply losing convergence guarantees or polynomial-time convergence.

This complexity gap naturally leads to the question of whether we can build FW-based algorithms with explicit guarantees on the active set cardinality beyond the trivial bound
$\min \{t, |\mathcal{V}|\}$ with $t$ the number of iterations and $\mathcal{V}$ the set of vertices of the feasible region (considered infinite if the region is not a polytope).\\

\citet{wirth2024pivotingframeworkfrankwolfealgorithms} introduced a general framework for active set size control in FW algorithms, ensuring that the number of vertices remains below the Carathéodory bound of $n+1$ vertices to form any iterate \citep{caratheodory1907variabilitatsbereich}, with $n$ the dimension of the feasible region which may be much lower than the ambient dimension, with an implementation building on top of \package{}.
The pivoting framework essentially maintains a sparse matrix representing active vertices, the weights of these vertices forming the current iterate $\vx_t$ being the solution to a linear system formed with this matrix and a right-hand side based on $\vx_t$.
The result compliments previous work by \citet{beck2017linearly} which also proposed an active set control, but further exploits potential sparsity of the iterates and offers more control on the accumulation of numerical errors.
\citet{wirth2024pivotingframeworkfrankwolfealgorithms} also established an optimal face identification property for BPCG on polytopes, as was shown for AFW by \citet{bomze2020active}.
Informally, the identification property states that if the objective function is convex and if the minimizer of the problem is unique and lies on a face of the polytope, then there exists a finite number of iterations after which all iterates lie on the face and do not leave it again.

\section{Conclusion}

This paper presents an overview of the last years of development of the \package{} package, with new algorithms developed and integrated since its inception, new oracles and step-size strategies.
We also highlighted a few applications of \package{}, both in the form of new packages leveraging it and research that leveraged it for computations.
The diversity of applications of both types shows the success of the package in its two targeted use cases: as a ready-to-use toolbox for computational research and as a platform for further algorithmic development of Frank-Wolfe-based methods.
Future work on the package and ecosystem will include further variants for generalized settings, and applications leveraging FW subproblems within a more complex task.

\section*{Acknowledgments}

The authors acknowledge the work of all contributors and the feedback from users of the package.
We thank in particular Victor Thouvenot for contributions to the Python wrapper and combinatorial oracle companion package, and Zev Woodstock for work on the block-coordinate Frank-Wolfe.
Research reported in this paper was partially supported through the Research Campus Modal funded by the German Federal Ministry of Education and Research (fund numbers 05M14ZAM,05M20ZBM) and the Deutsche Forschungsgemeinschaft (DFG) through the DFG Cluster of Excellence MATH+.
The work of Mathieu Besançon benefited from the support of the FMJH Program PGMO and from MIAI at Université Grenoble Alpes (grant ANR-19-P3IA-0003).
Part of this work was conducted while Sebastian Pokutta was visiting Tokyo University via a JSPS International Research Fellowship.

\bibliographystyle{icml2021}
\bibliography{refs}

\appendix

\section{Pseudo-code for the new Frank-Wolfe variants}\label{sec:pseudocode}

We present in this section the pseudo-code for the Frank-Wolfe variants introduced in the package.

\begin{algorithm}[ht]
    \caption{Blended Decomposition-Invariant Conditional Gradients}\label{alg:bdicg}
    \begin{algorithmic}[1]
        \Require Starting extreme point $\vx_0$.
        \For{$t \in 0\dots T$}
        \State $\vvv_t \gets \argmin_{\vvv \in \Xset} \innp{\nabla f(\vx_t)}{\vvv}$
        \State $\vs_t \gets \argmin_{\vvv \in \Fcal{\vx}} \innp{\nabla f(\vx_t)}{\vvv}$
        \State $\va_t \gets \argmax_{\vvv \in \Fcal{\vx}} \innp{\nabla f(\vx_t)}{\vvv}$
        \If{$\kappa \innp{\nabla f(\vx_t)}{\va_t - \vs_t} \geq \innp{\nabla f(\vx_t)}{\vx_t - \vvv_t}$} \Comment{perform in-face step}
        \State $\vd_t \gets \va_t - \vs_t$ \label{line:bdicg_criterion}
    \Else \Comment{global FW step}
        \State $\vd_t \gets \vx_t - \vvv_t$
        \EndIf
        \State $\gamma_{t} \gets \argmin_{\gamma \in [0,\gamma^{\max}_t]} f(\vx_t - \gamma_t \vd_t)$
        \State $\vx_{t+1} \gets \vx_t - \gamma_t \vd_t$
        \EndFor
        \State \Return $\vx_t$
    \end{algorithmic}
\end{algorithm}

In Line~\eqref{line:bdicg_criterion} of Algorithm~\ref{alg:bdicg}, we include a scaling factor $\kappa$ for the test that allows to promote sparsity, as already found atoms are reused, similar to Line~\eqref{bpcg:select} in Algorithm~\ref{alg:BPCG}. The default value is $\kappa = 2.0$ in the package.

\begin{algorithm}[ht]
    \caption{Blended Pairwise Conditional Gradients (BPCG)}\label{alg:BPCG}
    \begin{algorithmic}[1]
        \Require Convex smooth function $f$, starting vertex $\vx_0 \in \mathrm{vert}(\Xset{})$.
        \Ensure Points $\vx_{1}, \ldots, \vx_T \in \Xset{}$.
        \State $\mathcal{S}_{0} \gets \{ \vx_0 \}$
        \For{$t \in 0\dots T-1$}
        \State $\va_t \gets \argmax_{\vvv \in \mathcal{S}_t} \innp{\nabla f(\vx_t)}{\vvv}$ \Comment{away vertex}
        \State $\vs_t \gets \argmin_{\vvv \in \mathcal{S}_t} \innp{\nabla f(\vx_t)}{\vvv}$ \Comment{local FW}
        \State $\vvv_t \gets \argmin_{\vvv \in V(P)} \innp{\nabla f(\vx_t)}{\vvv}$ \Comment{global FW}
        \If{$\kappa \innp{\nabla f(\vx_t)}{\va_t - \vs_t} \geq \innp{\nabla f(\vx_t)}{\vx_t - \vvv_t}$} \Comment{pairwise step} \label{bpcg:select}
        \State $\vd_t \gets \va_t - \vs_t$
        \State $\lambda_{t}^{\ast} \gets c[\vx_{t}](\va_{t})$
        \State $\gamma_t \gets \argmin_{\gamma \in [0, \lambda_{t}^{\ast}]} f(\vx_{t} - \gamma \vd_{t})$
        \If{$\gamma_{t} < \lambda_{t}^{\ast}$} \Comment{descent step}
        \State $\mathcal{S}_{t+1} \gets \mathcal{S}_{t}$
    \Else \Comment{drop step}
        \State $\mathcal{S}_{t+1} \gets \mathcal{S}_{t} \setminus \{ \va_{t} \}$
        \EndIf
    \Else \Comment{FW step}
        \State $\vd_t \gets \vx_t - \vvv_t$
        \State $\gamma_t \gets \argmin_{\gamma \in [0, 1]} f(\vx_t - \gamma \vd_t)$
        \State $\mathcal{S}_{t+1} \gets \mathcal{S}_t \cup \{ \vvv_t \}$ \Comment{or $\mathcal{S}_{t+1} \gets \{\vvv_t\}$ if $\gamma_t = 1$}
        \EndIf
        \State $\vx_{t+1} \gets \vx_t - \gamma_t \vd_t$
        \EndFor
        \State \Return $\vx_T$
    \end{algorithmic}
\end{algorithm}

\section{Overview of use cases and connected ecosystem}\label{sec:application}

We next review a few use cases which leveraged \package{} for different applications, from structured machine learning to mixed-integer convex optimization.
We then highlight key elements of the ecosystem surrounding \package{}, namely the Python interface and the companion package for linear oracles.

\paragraph*{Interpretable Predictions with Rate-Distortion Explanations}

One of the prominent features of FW algorithms is the sparsity they induce in the iterates, in the sense of iterates being formed as a convex combination of a typically low number of extreme points of the feasible region.
When extreme points are themselves sparse, this directly results in sparse iterates in the traditional sense of a low support.
This sparsity of iterates not only holds at the optimum but throughout the whole solving process, making the method particularly suitable to large-scale learning problems for which the optimization procedure is typically terminated fairly early.
In \citet{macdonald2022interpretable}, this property is exploited within the \emph{Rate-Distortion Explanation} (RDE) framework.
Informally, this framework tackles the problem of finding a sparse set of relevant input features with which the prediction model can keep the expected training loss low,
while the irrelevant features are replaced with Gaussian noise. Denoting by $\vx$ the filter applied on features to quantify their relevance (zero being an irrelevant feature and one an important feature),
$D(\vx)$ is the distortion obtained when applying the filter $\vx$, with $D(\mathbf{1}) = 0$ by convention, and $D(\mathbf{0})$ being the expected loss of a (Gaussian) random input.
The RDE problem is expressed as:
\begin{align*}
    \text{($k$-RDE)}\;\; \min_{\vx}\; & D(\vx) \\
    \text{s.t. }\;                    & \vx \in \Xset = \{\vx \in \left[0,1\right]^n, \sum_i x_i \leq k \}.
\end{align*}
\Xset is the intersection of a $k$-sparse polytope with the nonnegative orthant, it is the convex hull of vertices formed as the sum of at most $k$ vectors of the standard basis.
The problem ($k$-RDE) requires a prior choice of $k$ controlling the sparsity level and does not provide a complete overview of the trade-off between sparsity of the explanation and predictive power.
A way to obtain an explanation for multiple sparsity levels is to rank features instead of simply activating or deactivating them.
Once features are ranked from the most to least relevant, a $k$-sparse explanation can be obtained by selecting the $k$ first features according to the ranking.
Optimizing an ordering of the feature is equivalent to finding a permutation matrix for the features in their initial order in $[n]$.
Optimizing over permutation matrices leads to hard combinatorial problems which are intractable at the scale considered for RDE.
We can relax the problem by optimizing over their convex hull, the Birkhoff polytope, consisting of all doubly stochastic matrices.
This results in the following problem:
\begin{align*}
    \text{(Ord-RDE)}\;\; \min_{\Pi \in B_n}\; & \frac{1}{n-1}  \sum_{k \in [n-1]}  D(\Pi \mathbf{p}_k),
\end{align*}
with $B_n$ the Birkhoff polytope for $n\times n$ matrices, $\mathbf{p}_k$ the vector with $k$ ones and $n-k$ zeros.
The problem (Ord-RDE) optimizes the ordering of features that will perform best on average over all rates.
Computational experiments from \citet{macdonald2022interpretable} showed that (Ord-RDE) performed well across all rates, close enough to ($k$-RDE)
trained for the specific rate and much better than an ordering obtained from the average of all rates consisting in $n-1$ distinct ($k$-RDE) problems.\\

A generic interface for Stochastic Frank-Wolfe (SFW) algorithms was in particular implemented to support this work, allowing users to define their own batch size and momentum,
with for instance the growing batch size algorithm of \citet{hazan2016variance} or the momentum of \citet{mokhtari2020stochastic}.
We refer the interested reader to \citet[Section 4.1]{braun2022conditional} for a recent review of SFW algorithms.

\paragraph*{\texttt{InferOpt.jl} and differentiable combinatorial layers}

Differentiable optimization has been a notorious recent advance in machine learning,
enabling the integration of arbitrary convex optimization problems as layers in neural networks,
see the seminal work in \citet{amos2017optnet,agrawal2019differentiating} for quadratic and conic layers and
\citet{agrawal2019differentiable,besanccon2024flexible,blondel2022efficient} for several implementations.
These approaches however fail to differentiate through linear optimization with respect to the objective function,
which can be expressed as:
\begin{align*}
    \argmin_{\vvv \in \Xset} \innp{\mathbf{\theta}}{\vvv},
\end{align*}
since the derivative of such problem w.r.t.~$\theta$ leads to a zero Jacobian almost everywhere when \Xset{} is a polytope.
An approach originally proposed in \citet{berthet2020learning,blondel2020learning} consists in regularizing
the problem with a random perturbation of $\theta$, resulting in a probability distribution over $\vvv \in \Xset$
or equivalently, to the addition of a strongly convex regularizer to the original combinatorial problem.
The work in \citet{dalle2022learning} applies FW methods to solve the regularized problem,
and derives the equivalent probability distribution over vertices from the weights of the vertices of \Xset{} in the final convex decomposition.
Their package \texttt{InferOpt.jl} leverages \package{} to solve the FW problems, extract the probability distribution and compute Jacobians of the regularized problems,
making combinatorial layers compatible with automatic differentiation to integrate them into learning pipelines.

\paragraph*{\texttt{Boscia.jl} and first-order mixed-integer convex optimization}

In \citet{hendrych2022convex}, \package{} was combined with a branch-and-bound framework to tackle a subclass of Mixed-Integer
Nonlinear Problems (MINLP)
\begin{align*}
    \min_{\vx\in\Xset} & \, f(\vx),
\end{align*}
where $f$ is a differentiable convex function and $\Xset$ is a nonconvex set containing linear, combinatorial, and integrality constraints.
The FW variants implemented in \package{} are used to solve the node problems with the distinction that the LMO solves the mixed-integer linear problem under additional bound constraints
\begin{align*}
    \min_{\vvv\in\Xset\cap[\vu,\vl]} & \, \innp{\nabla f(\vx)}{\vvv}.
\end{align*}
This results in obtaining integral feasible solutions, and thus a primal bound, from the root node.
Due to the expensiveness of the LMO, the lazification capabilities of the FW variants are heavily utilized.
The active set decomposition facilitates the warm-starting of the children nodes and the error adaptiveness enables a tightening of the
solving precision with an increasing depth in the tree.
Through the callback mechanism, the evaluation of a node can be stopped dynamically if either the dual bound has reached the primal bound of the tree or
there is a user-specified number of better candidate nodes to investigate.

\texttt{Boscia.jl} has successfully been applied to the Optimal Experiment Design Problem (OEDP) under the D-criterion and A-criterion in \citet{design_of_experiments_boscia_23}.
OEDP amounts to choosing a subset of experiments maximizing the information gain and under the two mentioned criteria leads to $\mathcal{NP}$-hard problems.

Another use case of \texttt{Boscia.jl} is the Network Design Problem for Traffic Assignment, see \citet{sharma2024network}.
The goal is to add edges to an existing design such that the traffic flow cost and the design cost are minimized.
On both of the aforementioned problems, \texttt{Boscia.jl} showed promising performance results.

\paragraph*{\texttt{BellPolytopes.jl} and large-scale projections for nonlocality}

In quantum information theory, more precisely in Bell nonlocality~\citep{brunner2014bell}, a question of interest is the membership problem for the so-called \emph{local polytope}.
In the bipartite setting, this polytope is affinely equivalent to the cut polytope~\citep{avis2006relationship} and practical instances are too large to be solved via facet enumeration or linear programming.
Instead, a FW approach can be taken to solve the approximate Carathéodory problem for this polytope~\citep{combettes2023revisiting}.
The corresponding works by \citet{designolle2023improved,designolle2024better,designolle2024symmetric} utilize BPCG and leading sparse decompositions and separating hyperplanes that beat the current state-of-the-art from the literature~\citep{hirsch2017better,divianszky2017qutrit}.
These previous results relied on a rediscovered standard FW algorithm~\citep{brierley2016convex} attributed to~\citet{gilbert1956iterative}.

Formally, the task at hand is to find the Euclidean distance between a tensor of size $m^N$ and the $N$-partite local correlation polytope which is defined by its vertices $d^{\vec{a}^{(1)}\ldots \vec{a}^{(N)}}$, which are rank-one tensors labelled by $\vec{a}^{(n)}=a^{(n)}_1\ldots a^{(n)}_m$ for $n\in[N]$ such that $a^{(n)}_{x_n}=\pm1$ and with entries expressed as:
\begin{align*}
    d^{\vec{a}^{(1)}\ldots \vec{a}^{(N)}}_{x_1\ldots x_N}=\prod_{n=1}^N a^{(n)}_{x_n}.
\end{align*}
Since the objective function is quadratic, its gradient is linear and the LMO reduces, in the bipartite case, to a Quadratic Unconstrained Binary Optimisation (QUBO) instance.
In practice, heuristic solutions for this LMO are good enough to provide progress throughout the algorithm.
On the one hand, when the focus is on the construction of an explicit decomposition of points lying inside of the local polytope, this approximate LMO suffices.
Note that the results in \citet{designolle2023improved} take a final extra step to convert the numerical output of the algorithm into an analytical decomposition of the point of interest at the expense of a slight geometrical contraction that can be of interest for similar membership problems.
On the other hand, for points outside of the local polytope, the LMO must be solved exactly once at the end of the algorithm.
In the multipartite case, symmetry could be exploited by \citet{designolle2024symmetric} to significantly reduce the dimension and accelerate the (exact) LMO by enumerating only the orbits of vertices.
Notice that the acceleration described in \cref{sec:quadratic} for quadratic objective functions was introduced to fit the needs of this use case.
We refer the reader interested in more hands-on details to the corresponding example in the documentation, entitled \emph{\href{https://zib-iol.github.io/FrankWolfe.jl/stable/examples/docs_12_quadratic_symmetric/}{Accelerations for quadratic functions and symmetric problems}}.

\paragraph*{\texttt{ApproximateVanishingIdeals.jl}}

In the context of classification, approximate vanishing ideal algorithms generate feature transformations into a high-dimensional space in which the dataset becomes linearly separable \citep{heldt2009approximate,limbeck2013computation,livni2013vanishing,wirth2022conditional,wirth2023approximate}.
The latter two works study the Oracle Approximate Vanishing Ideal (OAVI) algorithm, in which the feature transformation is created by solving several Least Absolute Shrinkage and Selection Operator (LASSO) problems with FW variants.
The solutions generated by the FW variants are sparse convex combinations of vertices of the feasible region, thus leading to robust and interpretable feature transformations.
The \texttt{ApproximateVanishingIdeals.jl} package \citep{kuzinowicz2023avi} implements most of the existing approximate vanishing ideal algorithms and relies on the \package{} implementations of the FW algorithms for solving the LASSO subproblems in OAVI.

\paragraph*{\texttt{KernelHerding.jl}}

The kernel herding algorithm is a method for learning in Markov random field models \citep{welling2009herding} and is equivalent to the FW algorithm \citep{bach2012equivalence}.
The only existing software option for kernel herding is the \texttt{KernelHerding.jl} package \citep{wirth2022kernelherding}, which builds on \package{}.

\paragraph*{Python interface}

We developed \href{https://github.com/ZIB-IOL/frankwolfe-py}{\texttt{frankwolfe-py}}, a Python interface exposing the same functionalities as the native Julia package.
All computations are run in Julia, offering the performance with little overhead.
The package lets users call the solver with \texttt{juliacall}, but also define their own custom oracles, objective function and gradient thanks to \texttt{PythonCall}.
We present below a simple example using the Python package.\\

\begin{lstlisting}[language=python]
from frankwolfepy import frankwolfe
import numpy as np

def f(x):
    return np.linalg.norm(x)**2

def grad(storage,x):
    for i in range(len(x)):
        storage[i] = x[i]

# Create the Linear Minimization Oracle
lmo_prob = frankwolfe.ProbabilitySimplexOracle(1)

# Compute first point
x0 = frankwolfe.compute_extreme_point(lmo_prob, np.zeros(5))

# Call the algorithm and get the result tuple
result = frankwolfe.frank_wolfe(f, grad, lmo_prob, x0)
\end{lstlisting}

\paragraph*{\texttt{CombinatorialLinearOracles.jl}}

As highlighted in \cref{sec:lmo}, many linear minimization oracles have been introduced.
Implementing some oracles directly in \package{} would however excessively increase the number of dependencies and load time, along with the scope of the package.
For these reasons, we introduced a companion package \texttt{CombinatorialLinearOracles.jl} to collect important combinatorial oracles, in particular for graph-related problems,
in which the corresponding oracle typically corresponds to finding a subset of edges or vertices of minimum or maximum weight.
Two examples already implemented are minimum-weight spanning trees using the Kruskal algorithm implemented in \texttt{Graphs.jl}
and minimum-weight matching using the \texttt{BlossomV} solver \citep{kolmogorov2009blossom} for the perfect matching subroutine and the reduction from \citet{matchings}.
A linear oracle for maximum-weight matching is of particular importance since it can be called in polynomial time in the size of the input graph, while any polytope representation has to be of exponential size \citep{rothvoss2017matching}.


\section{Benchmarks and computational evaluation}\label{sec:benchmarks}

Lastly, we have developed a benchmark repository\footnote{The benchmark repository is available at \url{https://github.com/ZIB-IOL/Benchmarks_for_FrankWolfe_and_Boscia}.} for \package{} using the \texttt{BenchmarkTools.jl} package.
We have benchmarked the original Frank-Wolfe algorithm (Standard), the lazified FW version (Lazy FW), Away-Frank-Wolfe (AFW), Pairwise Conditional Gradient (PCG), Blended Conditional Gradient (BCG), Blended Pairwise Conditional Gradient (BPCG), Decomposition Invariant Conditional Gradient (DICG) and Blended Decomposition Invariant Conditional Gradient (BDICG).
Note that the latter two variants can only be applied to problems over 0/1 polytopes and thus, will only be benchmarked on some of the problems.
Additionally, we have benchmarked the lazified versions of Standard, AFW, PCG, BCG and BPCG and the \texttt{ActiveSetQuadraticProductCaching} on the applicable problems.
The dual gap tolerance is set to the default value of \texttt{1e-7}.
The line search is also the default one, namely adaptive step size introduced in \cref{subsubsec:Adaptive}.
The benchmarks will be continuously extended to contain more problems and new variants.
The tables \cref{tab:A-Opt,tab:ActiveSetQuadraticProductCaching,tab:Birkhoff,tab:D-Opt,tab:Lazified,tab:Nuclear,tab:Poisson,tab:Simplex,tab:Sparse,tab:Spectrahedron} of the benchmarks in \cref{sec:benchmark_tables} display the geometric mean of the solving time and final dual gap over multiple samples and average estimated memory consumption.
Note that the memory values may be misleading because progress and hence allocations vastly differ between the variants.

In the following, we introduce the problems and showcase the performance of the different FW variants on them.

\paragraph{Birkhoff polytope}

Given a dimension $n\in\mathbb{N}$, we generate a fixed random matrix $\Tilde{\mathbf{X}}\in(0,1)^{n \times n}$ and consider the problem:
\begin{align*}
 \min_{\mathbf{X}} & \; \frac{1}{n^2}\|\mathbf{X} - \Tilde{\mathbf{X}}\|_F^2 \\
 \text{s.t.} & \; \mathbf{X} \in B_n,
\end{align*}
where $B_n$ denotes the Birkhoff polytope, i.e., the set of doubly stochastic matrices.
A call to the LMO corresponds to a call to the Hungarian algorithm computing an optimal linear assignment.

\paragraph{Nuclear norm ball for collaborative filtering}

For dimension $n\in\mathbb{N}$ and a given $k\in\mathbb{N}$, we generate a fixed matrix $\Tilde{\mathbf{X}} \in \mathbb{R}^{n\times n}$ of rank $k$.
After that, we randomly choose missing entries and consider the objective function:
\begin{align*}
    \min_{\mathbf{X}} & \; \frac{1}{2} \sum_{(i, j) \in P} (X_{ij} - \Tilde{X}_{ij})^2 \\
    \text{s.t.} & \; \|\mathbf{X}\|_{\mathrm{nuc}} \leq \tau,
\end{align*}
where $P$ is the set of present entries, $\|\cdot\|_{\mathrm{nuc}}$ is the nuclear norm, and $\tau$ is a hyperparameter.
The LMO for the nuclear norm ball consists in computing the negative-most singular vector pair of the gradient matrix.

\paragraph{The Optimal Experiment Design under the A-criterion and D-criterion}

The aim of the Optimal Experiment Design problem is to select the subset of experiments maximizing the information gain about the observed system.
The size of the subset is fixed and is usually significantly smaller than the number of total experiments.
The natural formulation is a MINLP but there is a continuous formulation called the limit problem whose feasible region is the standard probability simplex \citep{design_of_experiments_boscia_23}.
There are different information measures. Two popular ones are the A-criterion
\begin{align*}
\min_{\vx} & \; \tr\left(\left(\mathbf{A}^\intercal \diag(\vx)\mathbf{A}\right)^{-1}\right) \\
\text{s.t.} & \; \vx \in \Delta_n
\end{align*}
and the D-criterion
\begin{align*}
    \min_{\vx} & \; \log\det\left(\mathbf{A}^\intercal \diag(\vx)\mathbf{A}\right) \\
    \text{s.t.} & \; \vx \in \Delta_n.
\end{align*}
The matrix $\mathbf{A}$ encodes the a priori information about the experiments.
The variables $\vx$ can be interpreted as a probability distribution over experiments.

\paragraph{Poisson regression}

The Poisson regression problem aims to model count data where we assume that our data points $y_i$ follow a Poisson distribution.
The matrix $\mathbf{X}$ encodes the coefficients for the linear estimation of the expected values of $y_i$.
Note that $\vx_i$ denotes the $i$-th row of $\mathbf{X}$.

\begin{align*}
    \min_{\vw, b} & \; \sum\limits_{i=1}^n \exp\left(\vw^\intercal \vx_i + b\right) - y_i \left(\vw^\intercal \vx_i + b\right) + \alpha \norm{\vw}^2 \\
    \text{s.t.} & \;  -N \leq w_i \leq N \; \forall i \in [n] \\
    & \; b \in [-N, N].
\end{align*}

\paragraph{Probability simplex}

For dimensions $(m, n) \in \mathbb{N}^2$ we create a fixed random, normally distributed matrix $A\in\mathbb{R}^{m \times n}$ and a random normally distributed vector $b\in\mathbb{R}^m$ and consider the minimization problem
\begin{align*}
    \min_{\vx} & \; \frac{1}{m} \|A\vx + b\|_2^2\\
    \text{s.t.} & \; \vx \in \Delta_n.
\end{align*}

\paragraph{$K$-sparse polytope}

Given a dimension $n\in\mathbb{N}$ and a $K < n$, we generate a fixed random vector $\Tilde{x}\in[100]^{n}$. Let $\mathbf{y}=\frac{1}{\|\Tilde{\mathbf{x}}\|_1}\Tilde{\mathbf{x}}$ and consider the optimization problem:
\begin{align*}
    \min_{\vx} & \; \|\vx - \vy\|_2^2 \\
    \text{s.t.} & \; \vx \in \mathcal{B}_1(K) \cap \mathcal{B}_{\infty}
\end{align*}
where $\mathcal{B}_1(K)$ denotes the $\ell_1$-norm ball of radius $K$ and $\mathcal{B}_{\infty}$ the unit $\ell_\infty$-norm ball.
The LMO consists in selecting the $K$ lowest entries of the gradient.

\paragraph{Spectrahedron}

Lastly, we have a quadratic over the spectraplex $\mathcal{S}_{\tau}$ introduced in \cref{sec:lmo}.
For a matrix $\mathbf{Y}\in\mathbb{R}^{n\times n}$ and observed entries $P\subset[n]\times[n]$, we consider the problem:
\begin{align*}
    \min_{\mathbf{X}} & \; \frac{1}{2} \sum_{(i, j) \in P} (X_{ij} - Y_{ij})^2 \\
    \text{s.t.} & \; \mathbf{X} \in \mathcal{S}_{\tau}
\end{align*}
where $\tau$ is between 1.0 and 8.0. Similarly to the nuclear norm ball, the LMO consists in computing the eigenvector associated with the most negative eigenvalue of the gradient, assuming w.l.o.g.~that the function is invariant under transposes (thus preserves symmetries).

\begin{figure}[t]
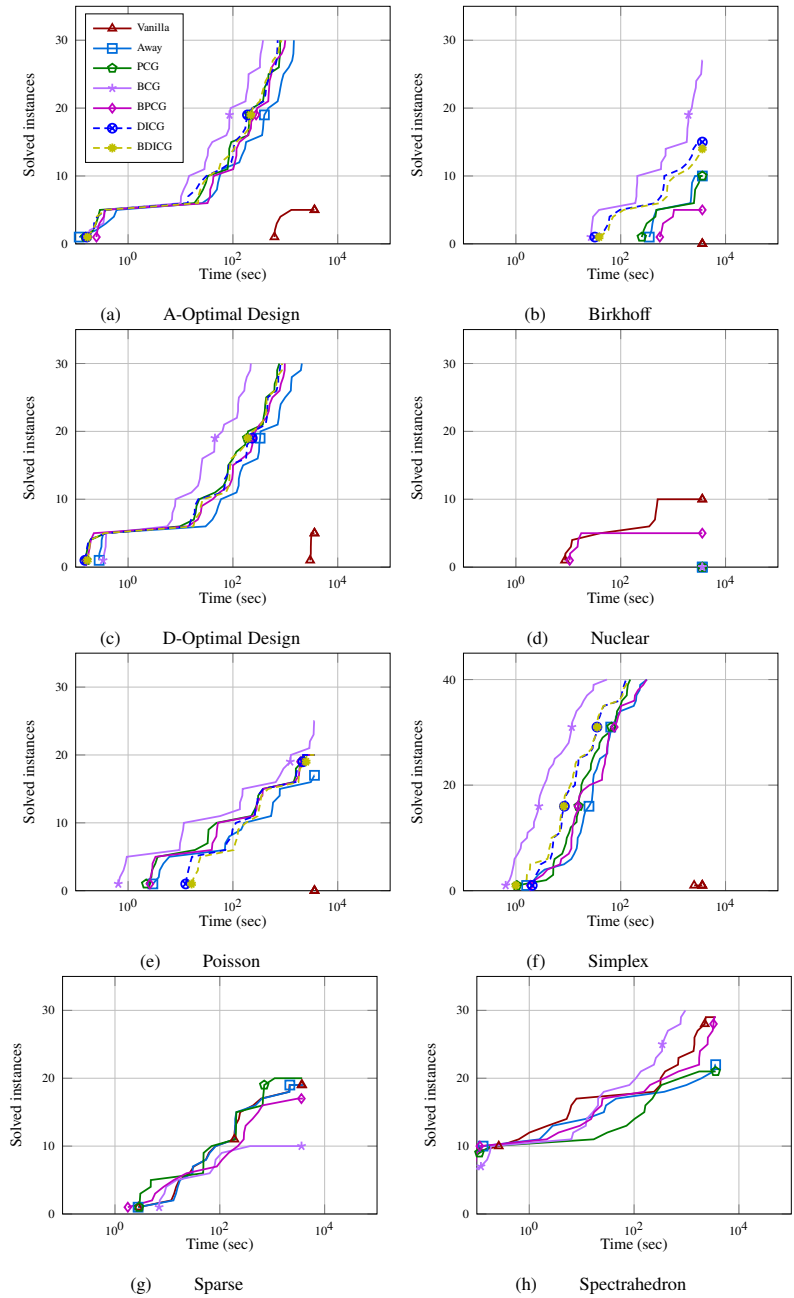

    \begin{center}
    \subfloat[][\hspace{0.5cm}A-Optimal Design]{
    \centering

    }
    \end{center}
    \caption{Number of solved instances for base variants.}
    \label{fig:base_solved_instances}
\end{figure}

\begin{figure}[ht]
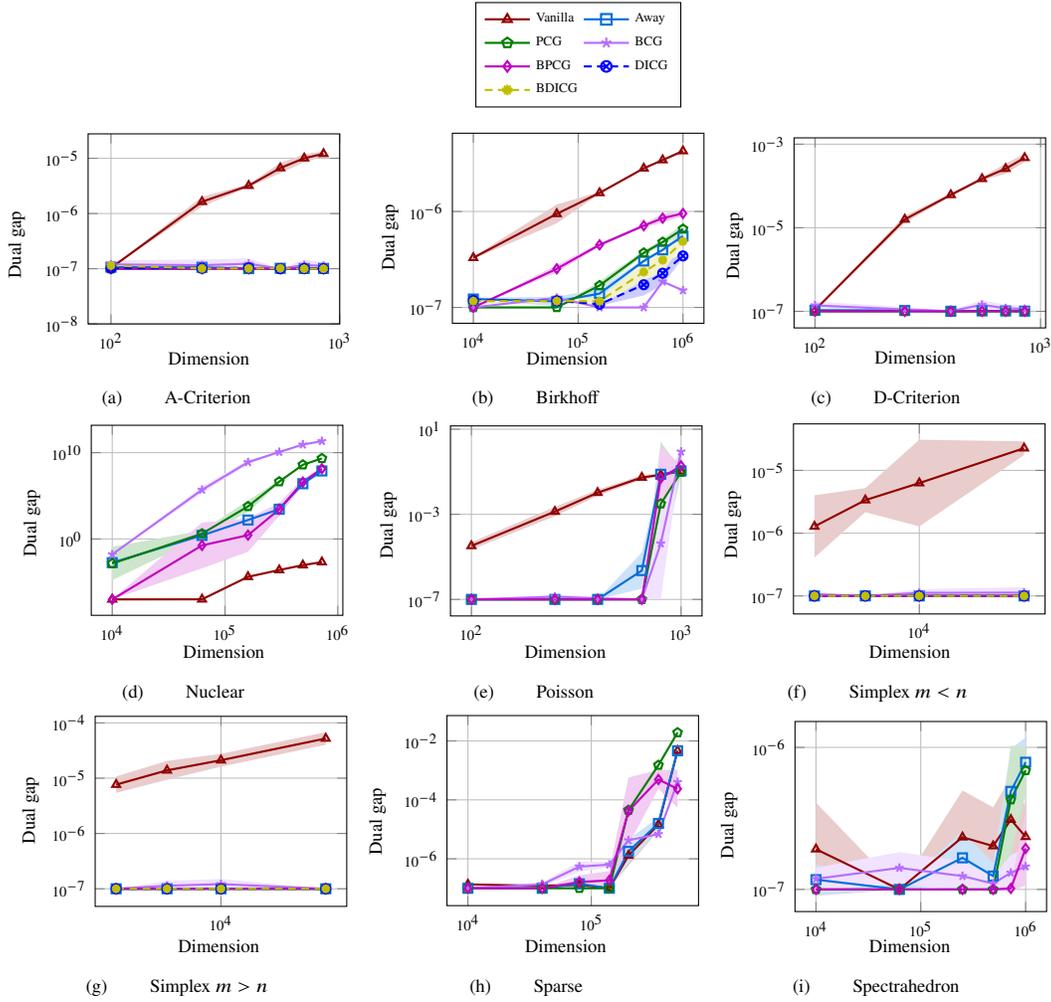

    \begin{center}
    \subfloat[][\hspace{0.5cm}A-Criterion]{
    \centering
    %
    }
    \end{center}
    \caption{Geometric mean of dual gaps of base variants.}
    \label{fig:base_dual_gaps}
\end{figure}


In \cref{fig:base_solved_instances}, we compare the number of optimally solved instances of the baseline variants for all problems.
The geometric mean and geometric standard deviation of the dual gaps are shown in \cref{fig:base_dual_gaps}.
The hardest problem is the Nuclear Norm Ball problem, see \cref{fig:base_solved_instances} (d), where most variants already timeout for the smaller dimension.
Standard FW shows the best performance, both in terms of solved instances and dual gaps.
This suggests that the active set is quite large and the FW vertex computed in any given iteration is not part of the active set at that point.

We observe a similar performance pattern for the Birkhoff problem, both Optimal Design problems, the Poisson problem and the Simplex problem.
While the decomposition-invariant methods outperform BPCG, especially for higher dimension, BCG is clearly the best performer.
This superior performance indicates that the weights of the vertices are shifted around considerably.
In this case, BCG has the advantages since it can move weights between multiple vertices simultaneously whereas both BPCG and the decomposition-invariant methods move weights only between two vertices.
Observe that on all variants, standard FW achieves small dual gap values even for the timed out instances on all five problems with exception of the Poisson instances of large dimension.

On the problem over the $K$-sparse polytope, standard FW and PCG showcase the best performance.
Note that the objective of the problem is strongly convex and the optimal solution is likely in the strict interior.
This together with the fact that the corresponding LMO is not very expensive explain the good performance of standard FW.
For the timed out instances, the quality of the achieved dual gap varies significantly.

There is no clear best performer on the Spectrahedron problem, see \cref{tab:Spectrahedron}.
Standard FW performs well for smaller dimensions.
The active set-based methods are superior for larger dimension, in particualr BPCG.
Considering that the he LMO of the Spectrahedron is comparatively expensive since it requires an eigenvalue decomposition, the observed results are to be expected.

\begin{figure}[ht]
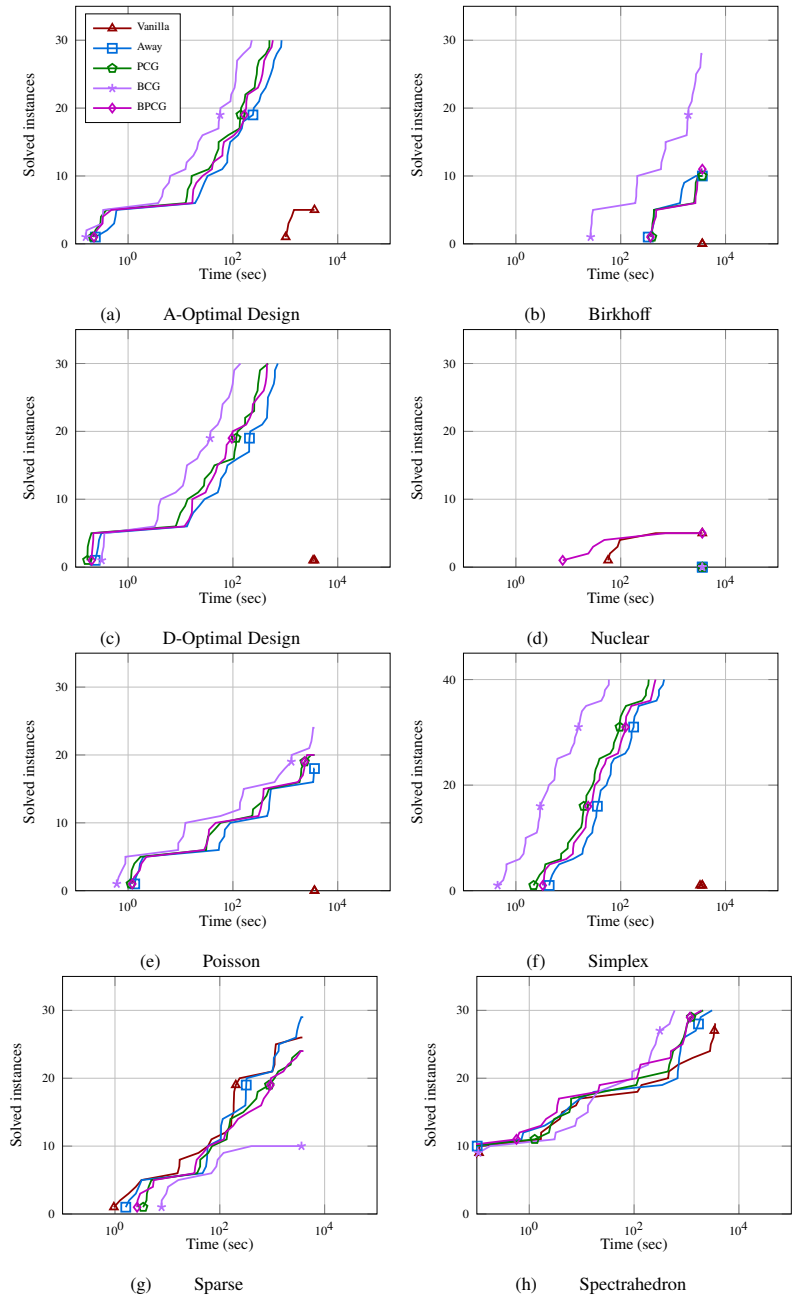

    \begin{center}
    \subfloat[][\hspace{0.5cm}A-Optimal Design]{
    \centering

    }
    \end{center}
    \caption{Number of solved instances for lazified variants.}
    \label{fig:lazy_solved_instances}

\end{figure}

\begin{figure}[ht]
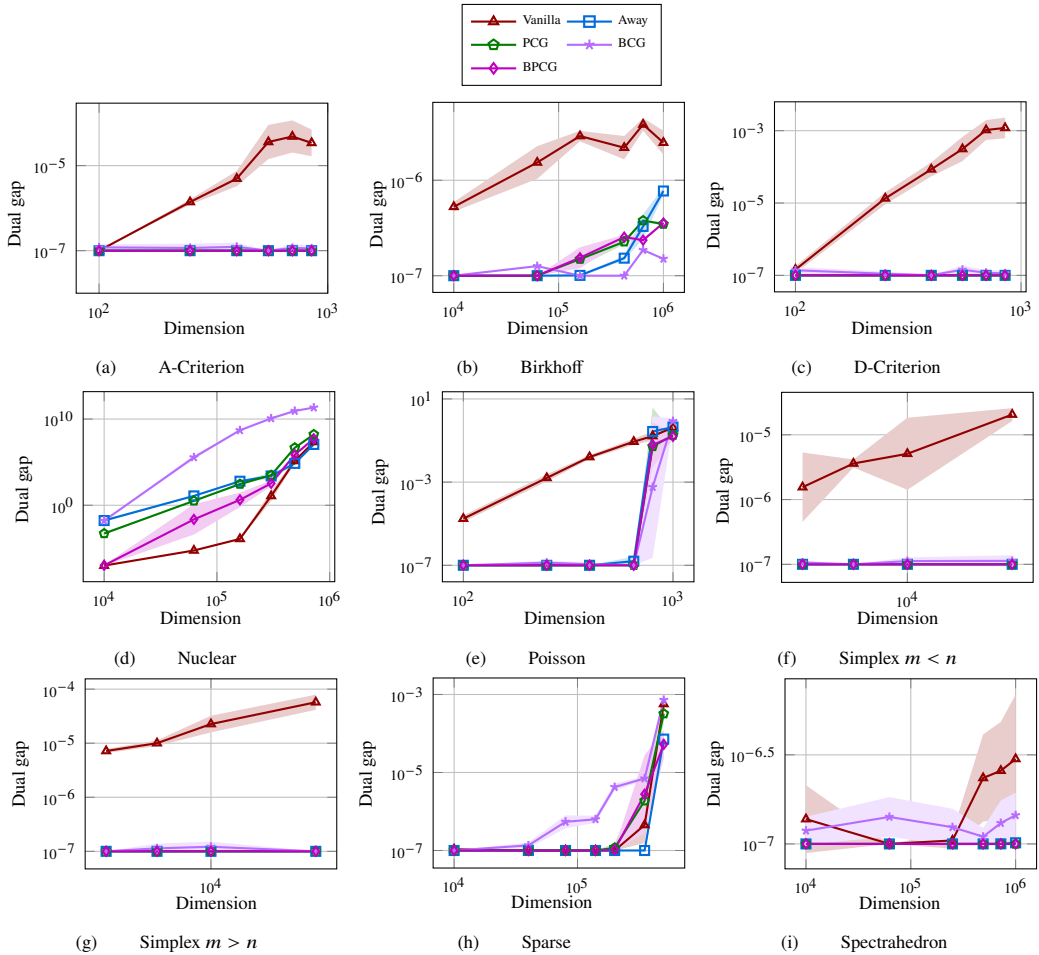

\begin{center}
\subfloat[][\hspace{0.5cm}A-Criterion]{
\centering
%
}
 \end{center}
\caption{Geometric mean and geometric standard deviation of dual gaps for lazified variants.}
\label{fig:lazy_dual_gaps}
\end{figure}

\cref{fig:lazy_solved_instances,fig:lazy_dual_gaps} showcase the performance of the lazified versions of the methods (there applicable).
We can cleary observe that lazification has a positive impact if the LMO of the corresponding problem is computionally expensive as is the case for the Spectrahedron, the Nuclear Norm Ball, the Poisson Problem and most noteably the $K$-Sparse Polytope.
It should be pointed out that the Poisson Problem was modelled using \texttt{MathOptInterface.jl}, hence the LMO performed by an LP solver instead of a specialised routine.

\begin{figure}[ht]
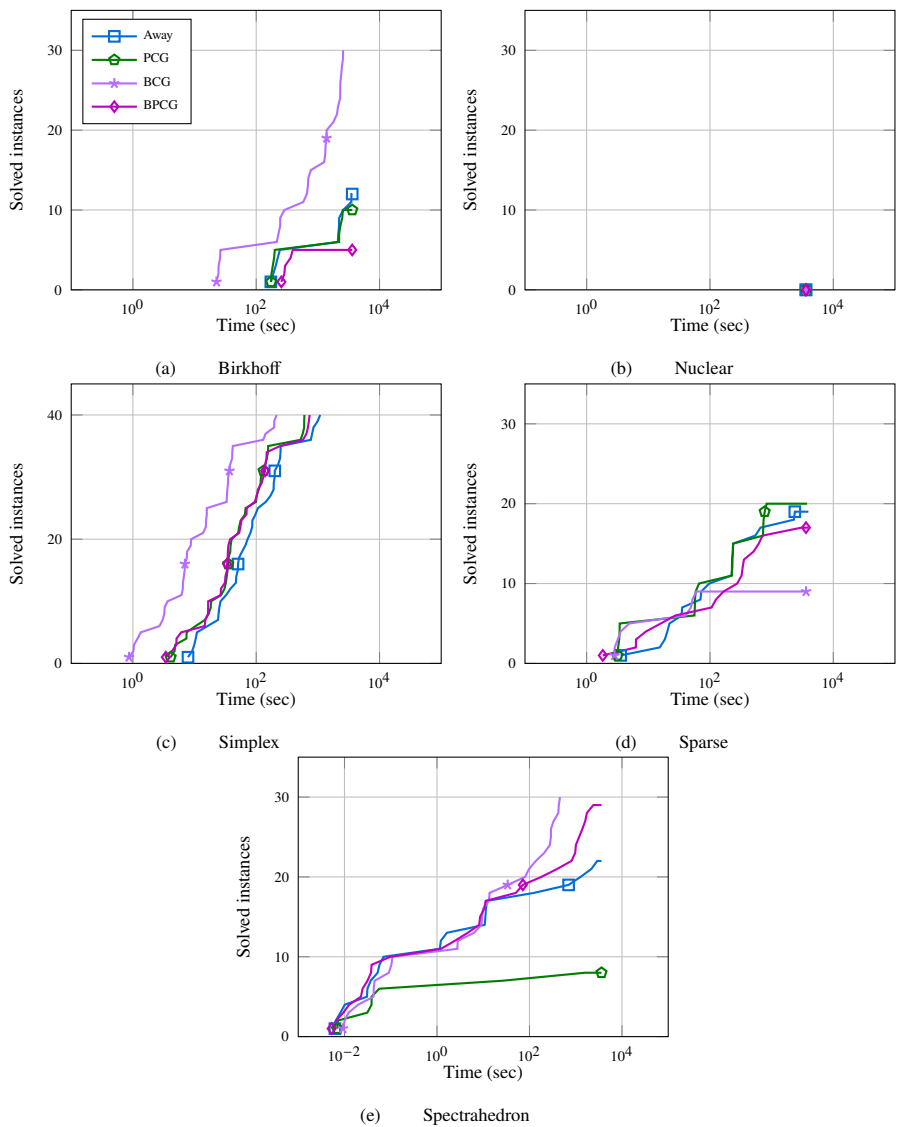

    \begin{center}
    \subfloat[][\hspace{0.5cm}Birkhoff]{
    \centering

    }
    \end{center}
    \caption{Number of solved instances for active set variants with product caching.}
    \label{fig:product_caching_solved_instances}
\end{figure}

\begin{figure}[ht]
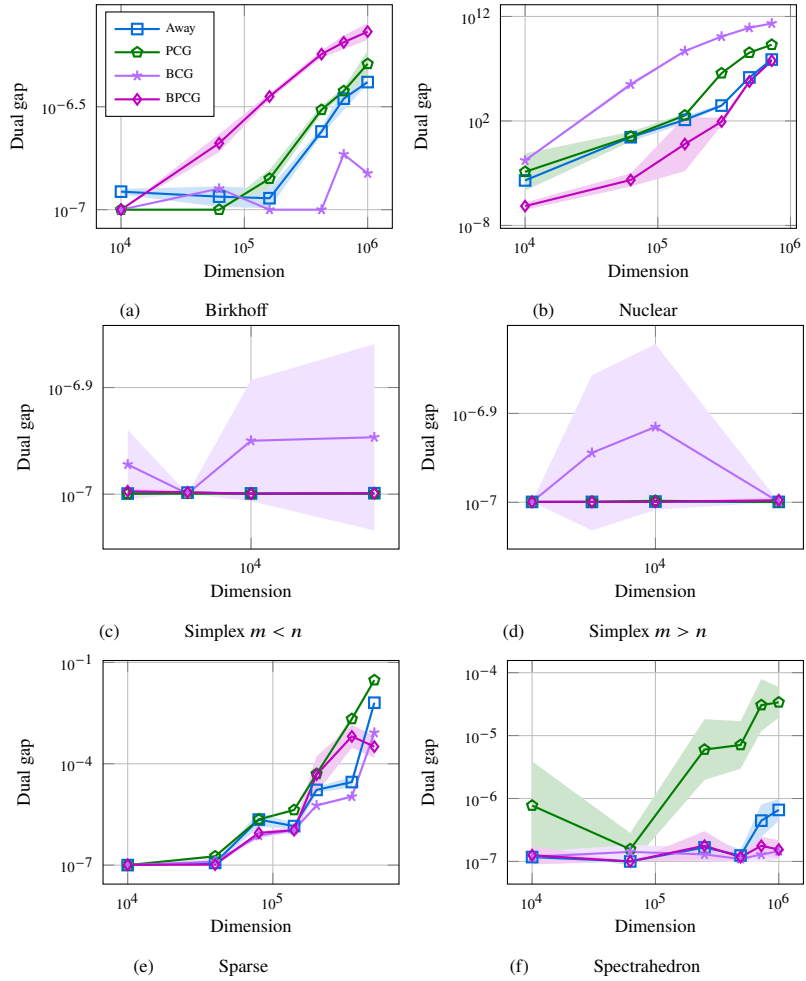

\begin{center}
\subfloat[][\hspace{0.8cm}Birkhoff]{
\centering
%
}
\end{center}
\caption{Geometric mean and geometric standard deviation of dual gaps for active set variants with product caching.}
\label{fig:product_caching_dual_gaps}
\end{figure}

In the benchmarks, the usage of \texttt{ActiveSetQuadraticProductCaching} is beneficial in cases where the Hessian of the corresponding objective function is trivial, i.e., the unit matrix.
This can be seen for the problem over the Spectrahedron, where the gain are quite substantial for BCG and BPCG \cref{fig:product_caching_solved_instances,fig:product_caching_dual_gaps}.
In contrast, the Simplex problem show a slow down, which can be understood when considering the cost of computing, as explained after \cref{eqn:argminmax}, many scalar products involving large matrix-vector multiplications: this introduces a significant overhead which only pays off after a certain number of iterations.
Note also that the number of atoms in the active set plays an important role in this trade-off: the more there are, the more interesting it can become to use \texttt{ActiveSetQuadraticProductCaching}.

To summarize, the choice of the most suitable variant depends on the cost of the LMO, the location of the optimal solution (strict interior or face), strong convexity of the objective, and geometry of the feasible region (polytope or nonlinear).
Standard FW is by design less memory-intensive but struggle if the LMO is more costly or the optimal solution is on a face.
The active set-based methods are the optimal choice for problems of larger dimensions.
While DICG requires less memory and provides directions well-aligned with the negative gradient, the second LMO call to compute the in-face vertex can cause computational overhead, as it is in many cases as expensive as the standard LMO.
Lazification can be a great addition if the LMO is very expensive and is, in particular, effective within BPCG.
Likewise, utilizing the \texttt{ActiveSetQuadraticProductCaching} can provide an advantage when optimizing a quadratic over the feasible region with expensive LMO.
The design of this diverse benchmark set precisely highlights the importance of the different FW methods for different types of problem classes, motivating their inclusion in the package.

\section{Benchmark tables}\label{sec:benchmark_tables}
\nopagebreak
\begin{sidewaystable}[ht]
    \centering
    \begin{adjustbox}{scale=0.7,center}
\begin{tabular}{c ccc ccc ccc ccc ccc ccc ccc}
            \toprule
            Dim & \multicolumn{3}{c}{Standard} & \multicolumn{3}{c}{AFW} & \multicolumn{3}{c}{PCG} & \multicolumn{3}{c}{BCG} & \multicolumn{3}{c}{BPCG} & \multicolumn{3}{c}{DICG} & \multicolumn{3}{c}{BDICG}\tabularnewline
            \cmidrule(lr){2-4}
            \cmidrule(lr){5-7}
            \cmidrule(lr){8-10}
            \cmidrule(lr){11-13}
            \cmidrule(lr){14-16}
            \cmidrule(lr){17-19}
            \cmidrule(lr){20-22}
            & \thead{Time \\ (s)} & \thead{Average \\ memory \\ (GB)} &\thead{Dual gap} &\thead{Time \\ (s)} & \thead{Average \\ memory \\ (GB)} &\thead{Dual gap} & \thead{Time \\ (s)} & \thead{Average \\ memory \\ (GB)} &\thead{Dual gap} & \thead{Time \\ (s)} & \thead{Average \\ memory \\ (GB)} &\thead{Dual gap} &\thead{Time \\ (s)} & \thead{Average \\ memory \\ (GB)}&\thead{Dual gap} & \thead{Time \\ (s)} & \thead{Average \\ memory \\ (GB)} &\thead{Dual gap} & \thead{Time \\ (s)} & \thead{Average \\ memory \\ (GB)} &\thead{Dual gap} \\
            \midrule
            100 & 7.76e2 & 1.30e3 & 1.03e-7 & 3.62e-1 & 6.93e-1 & 1.07e-7 & 2.32e-1 & 4.55e-1 & 1.00e-7 & 2.75e-1 & 2.64e-1 & 1.20e-7 & 2.97e-1 & 6.19e-1 & 1.00e-7 & 2.47e-1 & 4.95e-1 & 1.03e-7 & 2.45e-1 & 5.03e-1 & 1.14e-7 \\
            250 & 3.60e3* & 1.36e3 & 1.64e-6 & 5.14e1 & 1.42e1 & 1.06e-7 & 2.58e1 & 7.42 & 1.00e-7 & 1.18e1 & 3.36 & 1.16e-7 & 3.79e1 & 1.01e1 & 1.00e-7 & 1.95e1 & 7.47 & 1.02e-7 & 2.56e1 & 7.72 & 1.01e-7 \\
            400 & 3.60e3* & 1.73e3 & 3.19e-6 & 1.09e2 & 4.97e1 & 1.00e-7 & 8.45e1 & 3.73e1 & 1.00e-7 & 3.28e1 & 1.57e1 & 1.23e-7 & 1.15e2 & 4.99e1 & 1.00e-7 & 9.07e1 & 3.86e1 & 1.04e-7 & 8.00e1 & 3.83e1 & 1.01e-7 \\
            550 & 3.60e3* & 1.92e3 & 6.67e-6 & 3.97e2 & 2.08e2 & 1.02e-7 & 2.09e2 & 1.08e2 & 1.00e-7 & 8.41e1 & 4.85e1 & 1.00e-7 & 2.38e2 & 1.39e2 & 1.00e-7 & 1.82e2 & 1.14e2 & 1.00e-7 & 2.04e2 & 1.11e2 & 1.02e-7 \\
\bottomrule
        \end{tabular}
    \end{adjustbox}
    \caption{Performance comparison on the A-Criterion optimal design problem. Entries marked with * timed out.}
    \label{tab:A-Opt}

    \bigskip

    \begin{adjustbox}{scale=0.7,center}
        \begin{tabular}{c ccc ccc ccc ccc ccc ccc ccc}
            \toprule
            Dim & \multicolumn{3}{c}{Standard} & \multicolumn{3}{c}{AFW} & \multicolumn{3}{c}{PCG} & \multicolumn{3}{c}{BCG} & \multicolumn{3}{c}{BPCG} & \multicolumn{3}{c}{DICG} & \multicolumn{3}{c}{BDICG}\tabularnewline
            \cmidrule(lr){2-4}
            \cmidrule(lr){5-7}
            \cmidrule(lr){8-10}
            \cmidrule(lr){11-13}
            \cmidrule(lr){14-16}
            \cmidrule(lr){17-19}
            \cmidrule(lr){20-22}
            & \thead{Time \\ (s)} & \thead{Average \\ memory \\ (GB)} &\thead{Dual gap} &\thead{Time \\ (s)} & \thead{Average \\ memory \\ (GB)} &\thead{Dual gap} & \thead{Time \\ (s)} & \thead{Average \\ memory \\ (GB)} &\thead{Dual gap} & \thead{Time \\ (s)} & \thead{Average \\ memory \\ (GB)} &\thead{Dual gap} &\thead{Time \\ (s)} & \thead{Average \\ memory \\ (GB)}&\thead{Dual gap} & \thead{Time \\ (s)} & \thead{Average \\ memory \\ (GB)} &\thead{Dual gap} & \thead{Time \\ (s)} & \thead{Average \\ memory \\ (GB)} &\thead{Dual gap} \\
            \midrule
            10000 & 3.60e3* & 7.25e2 & 3.29e-7 & 4.09e2 & 1.33e2 & 1.22e-7 & 3.09e1 & 1.49e1 & 1.00e-7 & 3.43e2 & 1.05e2 & 1.00e-7 & 7.41e2 & 2.08e2 & 1.00e-7 & 5.59e1 & 1.47e1 & 1.18e-7 & 6.48e1 & 1.43e1 & 1.16e-7 \\
            62500 & 3.60e3* & 6.62e2 & 9.44e-7 & 2.30e3 & 4.89e2 & 1.16e-7 & 1.99e2 & 8.18e1 & 1.26e-7 & 2.80e3 & 5.48e2 & 1.00e-7 & 3.60e3* & 8.05e2 & 2.53e-7 & 6.04e2 & 8.12e1 & 1.17e-7 & 7.78e2 & 7.91e1 & 1.17e-7 \\
            160000 & 3.60e3* & 6.37e2 & 1.56e-6 & 3.60e3* & 5.23e2 & 1.39e-7 & 6.51e2 & 2.19e2 & 1.00e-7 & 3.60e3* & 5.42e2 & 1.69e-7 & 3.60e3* & 6.05e2 & 4.49e-7 & 2.14e3 & 1.77e2 & 1.08e-7 & 2.82e3 & 1.69e2 & 1.16e-7 \\
            422500 & 3.60e3* & 6.94e2 & 2.81e-6 & 3.60e3* & 4.50e2 & 3.05e-7 & 2.58e3 & 6.22e2 & 1.00e-7 & 3.60e3* & 4.74e2 & 3.70e-7 & 3.60e3* & 6.00e2 & 7.11e-7 & 3.60e3* & 2.44e2 & 1.73e-7 & 3.60e3* & 1.90e2 & 2.34e-7 \\
\bottomrule
        \end{tabular}
    \end{adjustbox}
    \caption{Performance comparison on the Birkhoff problem. Entries marked with * timed out.}
    \label{tab:Birkhoff}

    \bigskip

    \begin{adjustbox}{scale=0.7,center}
        \begin{tabular}{c ccc ccc ccc ccc ccc ccc ccc}
            \toprule
            Dim & \multicolumn{3}{c}{Standard} & \multicolumn{3}{c}{AFW} & \multicolumn{3}{c}{PCG} & \multicolumn{3}{c}{BCG} & \multicolumn{3}{c}{BPCG} & \multicolumn{3}{c}{DICG} & \multicolumn{3}{c}{BDICG}\tabularnewline
            \cmidrule(lr){2-4}
            \cmidrule(lr){5-7}
            \cmidrule(lr){8-10}
            \cmidrule(lr){11-13}
            \cmidrule(lr){14-16}
            \cmidrule(lr){17-19}
            \cmidrule(lr){20-22}
            & \thead{Time \\ (s)} & \thead{Average \\ memory \\ (GB)} &\thead{Dual gap} &\thead{Time \\ (s)} & \thead{Average \\ memory \\ (GB)} &\thead{Dual gap} & \thead{Time \\ (s)} & \thead{Average \\ memory \\ (GB)} &\thead{Dual gap} & \thead{Time \\ (s)} & \thead{Average \\ memory \\ (GB)} &\thead{Dual gap} &\thead{Time \\ (s)} & \thead{Average \\ memory \\ (GB)}&\thead{Dual gap} & \thead{Time \\ (s)} & \thead{Average \\ memory \\ (GB)} &\thead{Dual gap} & \thead{Time \\ (s)} & \thead{Average \\ memory \\ (GB)} &\thead{Dual gap} \\
            \midrule
            100 & 3.07e3 & 6.35e3 & 1.07e-7 & 3.03e-1 & 4.67e-1 & 1.06e-7 & 2.06e-1 & 3.18e-1 & 1.00e-7 & 3.62e-1 & 1.77e-1 & 1.38e-7 & 1.86e-1 & 4.03e-1 & 1.00e-7 & 2.02e-1 & 3.09e-1 & 1.06e-7 & 2.11e-1 & 3.39e-1 & 1.01e-7 \\
            250 & 3.60e3* & 1.44e3 & 1.58e-5 & 4.49e1 & 1.13e1 & 1.05e-7 & 1.73e1 & 5.96 & 1.00e-7 & 7.00 & 2.14 & 1.11e-7 & 2.37e1 & 7.19 & 1.00e-7 & 1.78e1 & 5.92 & 1.06e-7 & 2.19e1 & 6.05 & 1.03e-7 \\
            400 & 3.60e3* & 1.73e3 & 6.11e-5 & 1.34e2 & 5.98e1 & 1.00e-7 & 6.84e1 & 3.21e1 & 1.00e-7 & 2.27e1 & 1.04e1 & 1.01e-7 & 8.56e1 & 3.94e1 & 1.00e-7 & 7.40e1 & 3.24e1 & 1.04e-7 & 8.45e1 & 3.27e1 & 1.02e-7 \\
            550 & 3.60e3* & 1.95e3 & 1.48e-4 & 3.19e2 & 1.84e2 & 1.02e-7 & 1.48e2 & 9.73e1 & 1.00e-7 & 4.14e1 & 3.25e1 & 1.43e-7 & 2.17e2 & 1.16e2 & 1.00e-7 & 2.03e2 & 1.02e2 & 1.01e-7 & 1.56e2 & 1.03e2 & 1.00e-7 \\
\bottomrule
        \end{tabular}
    \end{adjustbox}
    \caption{Performance comparison on the D-Criterion optimal design problem. Entries marked with * timed out.}
    \label{tab:D-Opt}

    \bigskip

    \begin{adjustbox}{scale=0.7,center}
        \begin{tabular}{c ccc ccc ccc ccc ccc}
            \toprule
            Dim & \multicolumn{3}{c}{Standard} & \multicolumn{3}{c}{AFW} & \multicolumn{3}{c}{PCG} & \multicolumn{3}{c}{BCG} & \multicolumn{3}{c}{BPCG} \tabularnewline
            \cmidrule(lr){2-4}
            \cmidrule(lr){5-7}
            \cmidrule(lr){8-10}
            \cmidrule(lr){11-13}
            \cmidrule(lr){14-16}
            & \thead{Time \\ (s)} & \thead{Average \\ memory \\ (GB)} &\thead{Dual gap} &\thead{Time \\ (s)} & \thead{Average \\ memory \\ (GB)} &\thead{Dual gap} & \thead{Time \\ (s)} & \thead{Average \\ memory \\ (GB)} &\thead{Dual gap} & \thead{Time \\ (s)} & \thead{Average \\ memory \\ (GB)} &\thead{Dual gap} &\thead{Time \\ (s)} & \thead{Average \\ memory \\ (GB)}&\thead{Dual gap} \\
            \midrule
            10000 & 1.34e1 & 8.05e-1 & 1.00e-7 & 3.60e3* & 1.26e1 & 1.80e-3 & 3.60e3* & 1.86e1 & 1.44e-3 & 3.60e3* & 1.20e4 & 1.64e-2 & 1.36e1 & 3.41e-1 & 1.00e-7 \\
            62500 & 4.49e2 & 2.78 & 1.02e-7 & 3.60e3* & 1.47 & 2.59 & 3.60e3* & 1.50 & 4.11 & 3.60e3* & 9.94e3 & 4.96e5 & 3.60e3* & 1.76 & 1.79e-1 \\
            160000 & 3.60e3* & 2.23e1 & 4.12e-5 & 3.60e3* & 1.69 & 1.58e2 & 3.60e3* & 1.37 & 5.93e3 & 3.60e3* & 6.91e3 & 7.76e8 & 3.60e3* & 2.08 & 2.78 \\
            302500 & 3.60e3* & 2.45e1 & 2.52e-4 & 3.60e3* & 1.49 & 2.84e3 & 3.60e3* & 1.52 & 4.39e6 & 3.60e3* & 7.39e3 & 1.19e10 & 3.60e3* & 2.23 & 2.74e3 \\
            490000 & 3.60e3* & 2.26e1 & 9.25e-4 & 3.60e3* & 1.41 & 2.52e6 & 3.60e3* & 1.48 & 3.96e8 & 3.60e3* & 5.72e3 & 8.63e10 & 3.60e3* & 2.31 & 3.92e6 \\
            722500 & 3.60e3* & 1.00e1 & 2.09e-3 & 3.60e3* & 1.45 & 7.63e7 & 3.60e3* & 1.48 & 2.14e9 & 3.60e3* & 6.26e3 & 2.15e11 & 3.60e3* & 2.42 & 1.39e8 \\
            \bottomrule
        \end{tabular}
    \end{adjustbox}
    \caption{Performance comparison on the Nuclear problem. Entries marked with * timed out.}
    \label{tab:Nuclear}
\end{sidewaystable}

\begin{sidewaystable}
    \centering
    \begin{adjustbox}{scale=0.7,center}
        \begin{tabular}{c ccc ccc ccc ccc ccc ccc ccc}
            \toprule
            Dim & \multicolumn{3}{c}{Standard} & \multicolumn{3}{c}{AFW} & \multicolumn{3}{c}{PCG} & \multicolumn{3}{c}{BCG} & \multicolumn{3}{c}{BPCG} & \multicolumn{3}{c}{DICG} & \multicolumn{3}{c}{BDICG}\tabularnewline
            \cmidrule(lr){2-4}
            \cmidrule(lr){5-7}
            \cmidrule(lr){8-10}
            \cmidrule(lr){11-13}
            \cmidrule(lr){14-16}
            \cmidrule(lr){17-19}
            \cmidrule(lr){20-22}
            & \thead{Time \\ (s)} & \thead{Average \\ memory \\ (GB)} &\thead{Dual gap} &\thead{Time \\ (s)} & \thead{Average \\ memory \\ (GB)} &\thead{Dual gap} & \thead{Time \\ (s)} & \thead{Average \\ memory \\ (GB)} &\thead{Dual gap} & \thead{Time \\ (s)} & \thead{Average \\ memory \\ (GB)} &\thead{Dual gap} &\thead{Time \\ (s)} & \thead{Average \\ memory \\ (GB)}&\thead{Dual gap} & \thead{Time \\ (s)} & \thead{Average \\ memory \\ (GB)} &\thead{Dual gap} & \thead{Time \\ (s)} & \thead{Average \\ memory \\ (GB)} &\thead{Dual gap} \\
            \midrule
            100 & 3.60e3* & 2.11e3 & 3.27e-5 & 4.23 & 1.74 & 1.00e-7 & 2.92 & 1.98 & 1.00e-7 & 7.93e-1 & 9.38e-1 & 1.00e-7 & 2.88 & 2.08 & 1.00e-7 & 1.43e1 & 2.20 & 1.42e-7 & 2.02e1 & 2.51 & 1.42e-7 \\
            250 & 3.60e3* & 2.05e3 & 1.34e-3 & 1.00e2 & 4.57e1 & 1.00e-7 & 3.23e1 & 2.54e1 & 1.00e-7 & 1.05e1 & 1.01e1 & 1.34e-7 & 4.57e1 & 3.69e1 & 1.00e-7 & 8.96e1 & 3.44e1 & 1.83e-7 & 1.26e2 & 3.68e1 & 1.83e-7 \\
            400 & 3.60e3* & 1.31e3 & 1.04e-2 & 6.41e2 & 2.48e2 & 1.00e-7 & 2.88e2 & 1.54e2 & 1.00e-7 & 1.20e2 & 7.23e1 & 1.12e-7 & 3.07e2 & 1.72e2 & 1.00e-7 & 2.97e2 & 1.25e2 & 1.96e-7 & 3.48e2 & 1.30e2 & 1.96e-7 \\
            650 & 3.60e3* & 1.10e3 & 5.36e-2 & 3.49e3 & 1.24e3 & 2.27e-6 & 1.78e3 & 7.64e2 & 1.00e-7 & 9.47e2 & 4.74e2 & 1.00e-7 & 1.86e3 & 8.41e2 & 1.00e-7 & 1.82e3 & 7.38e2 & 1.58e-7 & 2.08e3 & 7.59e2 & 1.58e-7 \\
\bottomrule
        \end{tabular}
    \end{adjustbox}
    \caption{Performance comparison on the Poisson problem. Entries marked with * timed out.}
    \label{tab:Poisson}

    \bigskip

    \begin{adjustbox}{scale=0.7,center}
        \begin{tabular}{c ccc ccc ccc ccc ccc ccc ccc}
            \toprule
            Dim & \multicolumn{3}{c}{Standard} & \multicolumn{3}{c}{AFW} & \multicolumn{3}{c}{PCG} & \multicolumn{3}{c}{BCG} & \multicolumn{3}{c}{BPCG} & \multicolumn{3}{c}{DICG} & \multicolumn{3}{c}{BDICG}\tabularnewline
            \cmidrule(lr){2-4}
            \cmidrule(lr){5-7}
            \cmidrule(lr){8-10}
            \cmidrule(lr){11-13}
            \cmidrule(lr){14-16}
            \cmidrule(lr){17-19}
            \cmidrule(lr){20-22}
            & \thead{Time \\ (s)} & \thead{Average \\ memory \\ (GB)} &\thead{Dual gap} &\thead{Time \\ (s)} & \thead{Average \\ memory \\ (GB)} &\thead{Dual gap} & \thead{Time \\ (s)} & \thead{Average \\ memory \\ (GB)} &\thead{Dual gap} & \thead{Time \\ (s)} & \thead{Average \\ memory \\ (GB)} &\thead{Dual gap} &\thead{Time \\ (s)} & \thead{Average \\ memory \\ (GB)}&\thead{Dual gap} & \thead{Time \\ (s)} & \thead{Average \\ memory \\ (GB)} &\thead{Dual gap} & \thead{Time \\ (s)} & \thead{Average \\ memory \\ (GB)} &\thead{Dual gap} \\
            \midrule
            5000, 2500 & 3.35e3 & 1.03e3 & 1.28e-6 & 3.18 & 7.44e-1 & 1.00e-7 & 3.55 & 4.18e-1 & 1.00e-7 & 8.53e-1 & 1.40e-1 & 1.07e-7 & 3.77 & 4.48e-1 & 1.00e-7 & 2.48 & 4.12e-1 & 1.00e-7 & 1.57 & 4.07e-1 & 1.00e-7 \\
            5000, 6000 & 3.60e3* & 3.06e2 & 7.64e-6 & 1.70e1 & 1.04 & 1.00e-7 & 1.66e1 & 5.95e-1 & 1.00e-7 & 2.56 & 2.10e-1 & 1.00e-7 & 1.36e1 & 6.51e-1 & 1.00e-7 & 7.53 & 5.84e-1 & 1.00e-7 & 7.23 & 5.76e-1 & 1.00e-7 \\
            7000, 3500 & 3.60e3* & 6.95e2 & 3.36e-6 & 1.69e1 & 1.13 & 1.00e-7 & 8.18 & 6.47e-1 & 1.00e-7 & 1.30 & 2.20e-1 & 1.00e-7 & 9.49 & 6.95e-1 & 1.00e-7 & 4.93 & 6.37e-1 & 1.00e-7 & 4.39 & 6.31e-1 & 1.00e-7 \\
            7000, 8500 & 3.60e3* & 4.96e2 & 1.38e-5 & 2.86e1 & 3.52 & 1.00e-7 & 2.44e1 & 1.95 & 1.00e-7 & 3.84 & 6.18e-1 & 1.14e-7 & 4.79e1 & 2.09 & 1.00e-7 & 1.44e1 & 1.94 & 1.00e-7 & 1.33e1 & 1.91 & 1.00e-7 \\
\bottomrule
        \end{tabular}
    \end{adjustbox}
    \caption{Performance comparison on the Simplex problem. Entries marked with * timed out.}
    \label{tab:Simplex}

    \bigskip

    \begin{adjustbox}{scale=0.7,center}
        \begin{tabular}{c ccc ccc ccc ccc ccc}
            \toprule
            Dim & \multicolumn{3}{c}{Standard} & \multicolumn{3}{c}{AFW} & \multicolumn{3}{c}{PCG} & \multicolumn{3}{c}{BCG} & \multicolumn{3}{c}{BPCG} \tabularnewline
            \cmidrule(lr){2-4}
            \cmidrule(lr){5-7}
            \cmidrule(lr){8-10}
            \cmidrule(lr){11-13}
            \cmidrule(lr){14-16}
            & \thead{Time \\ (s)} & \thead{Average \\ memory \\ (GB)} &\thead{Dual gap} &\thead{Time \\ (s)} & \thead{Average \\ memory \\ (GB)} &\thead{Dual gap} & \thead{Time \\ (s)} & \thead{Average \\ memory \\ (GB)} &\thead{Dual gap} & \thead{Time \\ (s)} & \thead{Average \\ memory \\ (GB)} &\thead{Dual gap} &\thead{Time \\ (s)} & \thead{Average \\ memory \\ (GB)}&\thead{Dual gap} \\
            \midrule
            10000 & 1.35e1 & 4.86e1 & 1.35e-7 & 1.40e1 & 4.88e1 & 1.00e-7 & 3.51 & 1.12e1 & 1.00e-7 & 9.06 & 1.07e1 & 1.00e-7 & 5.82 & 2.12e1 & 1.00e-7 \\
            40000 & 5.52e1 & 1.97e2 & 1.20e-7 & 5.66e1 & 2.00e2 & 1.00e-7 & 5.14e1 & 1.75e2 & 1.00e-7 & 1.11e2 & 1.59e2 & 1.35e-7 & 1.16e2 & 3.95e2 & 1.00e-7 \\
            80000 & 9.09e2 & 3.03e3 & 1.30e-7 & 9.23e2 & 3.07e3 & 1.42e-7 & 2.01e2 & 7.44e2 & 1.00e-7 & 3.60e3* & 1.21e4 & 5.44e-7 & 1.30e3 & 4.36e3 & 1.63e-7 \\
            140000 & 2.10e2 & 7.37e2 & 1.00e-7 & 1.92e2 & 7.56e2 & 1.00e-7 & 7.40e2 & 2.28e3 & 1.00e-7 & 3.60e3* & 1.12e4 & 6.34e-7 & 5.44e2 & 1.99e3 & 1.86e-7 \\
            200000 & 3.60e3* & 1.26e4 & 1.29e-6 & 3.60e3* & 1.21e4 & 1.78e-6 & 3.60e3* & 1.17e4 & 4.46e-5 & 3.60e3* & 1.26e4 & 4.25e-6 & 3.60e3* & 1.22e4 & 4.16e-5 \\
            350000 & 3.60e3* & 1.19e4 & 1.46e-5 & 3.60e3* & 1.19e4 & 1.58e-5 & 3.60e3* & 1.12e4 & 1.51e-3 & 3.60e3* & 1.15e4 & 6.90e-6 & 3.60e3* & 1.20e4 & 4.92e-4 \\
            500000 & 3.60e3* & 1.14e4 & 4.64e-3 & 3.60e3* & 1.17e4 & 4.63e-3 & 3.60e3* & 1.11e4 & 1.92e-2 & 3.60e3* & 1.10e4 & 4.01e-4 & 3.60e3* & 1.16e4 & 2.36e-4 \\
            \bottomrule
        \end{tabular}
    \end{adjustbox}
    \caption{Performance comparison on the Sparse problem. Entries marked with * timed out.}
    \label{tab:Sparse}

    \bigskip

    \begin{adjustbox}{scale=0.7,center}
        \begin{tabular}{c ccc ccc ccc ccc ccc}
            \toprule
            Dim & \multicolumn{3}{c}{Standard} & \multicolumn{3}{c}{AFW} & \multicolumn{3}{c}{PCG} & \multicolumn{3}{c}{BCG} & \multicolumn{3}{c}{BPCG} \tabularnewline
            \cmidrule(lr){2-4}
            \cmidrule(lr){5-7}
            \cmidrule(lr){8-10}
            \cmidrule(lr){11-13}
            \cmidrule(lr){14-16}
            & \thead{Time \\ (s)} & \thead{Average \\ memory \\ (GB)} &\thead{Dual gap} &\thead{Time \\ (s)} & \thead{Average \\ memory \\ (GB)} &\thead{Dual gap} & \thead{Time \\ (s)} & \thead{Average \\ memory \\ (GB)} &\thead{Dual gap} & \thead{Time \\ (s)} & \thead{Average \\ memory \\ (GB)} &\thead{Dual gap} &\thead{Time \\ (s)} & \thead{Average \\ memory \\ (GB)}&\thead{Dual gap} \\
            \midrule
            10000 & 5.05e-2 & 6.41e-3 & 1.92e-7 & 1.52e-2 & 2.85e-3 & 1.17e-7 & 1.53e-2 & 2.88e-3 & 1.00e-7 & 1.55e-2 & 2.33e-2 & 1.19e-7 & 2.17e-2 & 3.51e-3 & 1.00e-7 \\
            62500 & 4.15e-2 & 1.05e-2 & 1.00e-7 & 6.81e-2 & 1.08e-2 & 1.00e-7 & 8.42e-2 & 1.11e-2 & 1.00e-7 & 1.04e-1 & 4.24e-2 & 1.42e-7 & 5.81e-2 & 1.07e-2 & 1.00e-7 \\
            250000 & 2.02e1 & 4.18e-1 & 2.33e-7 & 3.27e1 & 2.27e-1 & 1.67e-7 & 7.83e1 & 1.44e-1 & 1.00e-7 & 1.63e1 & 6.33 & 1.24e-7 & 2.29e1 & 1.87e-1 & 1.00e-7 \\
            490000 & 1.28e1 & 3.02e-1 & 2.02e-7 & 4.55e1 & 3.31e-1 & 1.24e-7 & 2.19e2 & 2.91e-1 & 1.00e-7 & 2.46e1 & 1.08e1 & 1.10e-7 & 3.18e1 & 3.17e-1 & 1.00e-7 \\
            722500 & 6.98e2 & 3.27 & 3.10e-7 & 3.52e3 & 2.39 & 4.90e-7 & 3.12e3 & 1.35 & 4.29e-7 & 2.71e2 & 4.70e2 & 1.31e-7 & 1.70e3 & 8.13 & 1.02e-7 \\
            1000000 & 1.93e3 & 7.63 & 2.35e-7 & 3.60e3* & 2.45 & 7.86e-7 & 3.60e3* & 1.53 & 6.88e-7 & 5.81e2 & 9.73e2 & 1.45e-7 & 2.94e3 & 1.34e1 & 1.94e-7 \\
            \bottomrule
        \end{tabular}
    \end{adjustbox}
    \caption{Performance comparison on the Spectrahedron problem. Entries marked with * timed out.}
    \label{tab:Spectrahedron}
\end{sidewaystable}

\begin{table}
    \centering
    \begin{adjustbox}{scale=0.60,center}
        \begin{tabular}{c c ccc ccc ccc ccc ccc }
            \toprule
            Problem & Dim & \multicolumn{3}{c}{Standard} & \multicolumn{3}{c}{Away} & \multicolumn{3}{c}{PCG} & \multicolumn{3}{c}{BCG} & \multicolumn{3}{c}{BPCG} \tabularnewline
            \cmidrule(lr){3-5}
            \cmidrule(lr){6-8}
            \cmidrule(lr){9-11}
            \cmidrule(lr){12-14}
            \cmidrule(lr){15-17}
            & & \thead{Time \\ (s)} & \thead{Average \\ memory \\ (GB)} &\thead{Dual gap} &\thead{Time \\ (s)} & \thead{Average \\ memory \\ (GB)} &\thead{Dual gap} & \thead{Time \\ (s)} & \thead{Average \\ memory \\ (GB)} &\thead{Dual gap} & \thead{Time \\ (s)} & \thead{Average \\ memory \\ (GB)} &\thead{Dual gap} &\thead{Time \\ (s)} & \thead{Average \\ memory \\ (GB)}&\thead{Dual gap} \\
            \midrule
            A-Criterion & $n$ & & & & & & & & & & & & & & &  \\
            \midrule
            & 100 & 1.20e3 & 2.29e3 & 1.00e-7 & 3.62e-1 & 9.16e-1 & 1.07e-7 & 2.32e-1 & 5.84e-1 & 1.00e-7 & 2.75e-1 & 2.64e-1 & 1.20e-7 & 2.97e-1 & 6.92e-1 & 1.00e-7 \\
            & 250 & 3.60e3* & 1.46e3 & 1.39e-6 & 5.14e1 & 1.48e1 & 1.06e-7 & 2.58e1 & 8.69 & 1.00e-7 & 1.18e1 & 3.36 & 1.16e-7 & 3.79e1 & 1.07e1 & 1.00e-7 \\
            & 400 & 3.60e3* & 1.58e3 & 4.95e-6 & 1.09e2 & 7.46e1 & 1.00e-7 & 8.45e1 & 4.30e1 & 1.00e-7 & 3.28e1 & 1.57e1 & 1.23e-7 & 1.15e2 & 4.93e1 & 1.00e-7 \\
            & 550 & 3.60e3* & 2.02e3 & 3.60e-5 & 3.97e2 & 2.20e2 & 1.02e-7 & 2.09e2 & 1.24e2 & 1.00e-7 & 8.41e1 & 4.85e1 & 1.00e-7 & 2.38e2 & 1.44e2 & 1.00e-7 \\
            & 700 & 3.60e3* & 2.16e3 & 4.84e-5 & 8.10e2 & 4.88e2 & 1.01e-7 & 4.20e2 & 2.80e2 & 1.00e-7 & 1.88e2 & 1.19e2 & 1.16e-7 & 5.08e2 & 3.21e2 & 1.00e-7 \\
\bottomrule
        \end{tabular}
    \end{adjustbox}
    \caption{Performance comparison of lazified variants. Entries marked with * timed out.}
    \label{tab:Lazified}
\end{table}

\begin{table}
    \centering
    \begin{adjustbox}{scale=0.70,center}
        \begin{tabular}{c c ccc ccc ccc ccc }
            \toprule
            Problem & Dim & \multicolumn{3}{c}{Away} & \multicolumn{3}{c}{PCG} & \multicolumn{3}{c}{BCG} & \multicolumn{3}{c}{BPCG} \tabularnewline
            \cmidrule(lr){3-5}
            \cmidrule(lr){6-8}
            \cmidrule(lr){9-11}
            \cmidrule(lr){12-14}
            & &\thead{Time \\ (s)} & \thead{Average \\ memory \\ (GB)} &\thead{Dual gap} & \thead{Time \\ (s)} & \thead{Average \\ memory \\ (GB)} &\thead{Dual gap} & \thead{Time \\ (s)} & \thead{Average \\ memory \\ (GB)} &\thead{Dual gap} &\thead{Time \\ (s)} & \thead{Average \\ memory \\ (GB)}&\thead{Dual gap} \\\\
            \midrule
            Birkhoff & $n$ & & & & & & & & & & & &  \\
            \midrule
            & 10000 & 2.05e2 & 9.96e1 & 1.22e-7 & 1.86e2 & 1.20e2 & 1.00e-7 & 2.48e1 & 3.44e1 & 1.00e-7 & 3.14e2 & 1.45e2 & 1.00e-7 \\
            & 62500 & 2.24e3 & 3.46e2 & 1.16e-7 & 2.35e3 & 5.74e2 & 1.00e-7 & 2.33e2 & 1.79e2 & 1.26e-7 & 3.60e3* & 6.46e2 & 2.10e-7 \\
            & 160000 & 3.56e3 & 4.58e2 & 1.14e-7 & 3.60e3* & 7.13e2 & 1.41e-7 & 6.75e2 & 5.32e2 & 1.00e-7 & 3.60e3* & 5.21e2 & 3.55e-7 \\
            & 422500 & 3.60e3* & 3.23e2 & 2.40e-7 & 3.60e3* & 5.19e2 & 3.06e-7 & 2.23e3 & 1.13e3 & 1.00e-7 & 3.60e3* & 3.93e2 & 5.69e-7 \\
\bottomrule
        \end{tabular}
    \end{adjustbox}
    \caption{Performance comparison of the active set variants with \texttt{ActiveSetQuadraticProductCaching}. Entries marked with * timed out.}
    \label{tab:ActiveSetQuadraticProductCaching}
\end{table}

\end{document}